\documentclass{amsart}

\usepackage{latexsym,amsthm,amssymb}

\usepackage{bbm}

\DeclareMathAlphabet{\mathpzc}{OT1}{pzc}{m}{it}

\usepackage[leqno]{amsmath}

\usepackage{mathrsfs}

\usepackage[all]{xy}

\usepackage{color}
\usepackage{colordvi}

\usepackage{verbatim}
\usepackage{ifthen}
\newboolean{hide.proofs}
\setboolean{hide.proofs}{false}
\newtheorem{theorem}{Theorem}[section]
\newtheorem{lemma}[theorem]{Lemma}

\newtheorem{proposition}[theorem]{Proposition}
\newtheorem{cor}[theorem]{Corollary}
\newtheorem{remark}[theorem]{Remark}
\newtheorem{noname}[theorem]{}

\newtheorem*{acknowledgement}{Acknowledgements}

\newtheorem{lemma-conjecture}[theorem]{Lemma--Conjecture}

\newtheorem{corollary}[theorem]{Corollary}

\newtheorem{question}[theorem]{Question}

\usepackage[hyperindex=true,plainpages=false,colorlinks=true,]{hyperref}
\numberwithin{equation}{theorem}

\renewcommand{\mathcal}{\mathscr}

\newcommand{\SB}{{\mathcal{B}}}

\newcommand{\SE}{{\mathcal{E}}}

\newcommand{\SL}{{\mathcal{L}}}

\newcommand{\SN}{{\mathcal{N}}}
\newcommand{\SO}{{\mathcal{O}}}

\newcommand{\ST}{{\mathcal{T}}}

\newcommand{\SX}{{\mathcal{X}}}
\newcommand{\SY}{{\mathcal{Y}}}

\newcommand{\PP}{\mathbf{P}}

\newcommand{\WA}{\widetilde{A}}

\newcommand{\WL}{\widetilde{L}}

\begin{document}

\title{Deformation of canonical morphisms and the moduli of surfaces of general type}

\author{Francisco Javier Gallego}
\author{Miguel Gonz\'alez}
\author{\\ Bangere P. Purnaprajna}

\address{Departamento de \'Algebra, Universidad Complutense de Madrid}
\email{gallego@mat.ucm.es}
\address{Departamento de \'Algebra, Universidad Complutense de Madrid}
\email{mgonza@mat.ucm.es}
\address{Department of Mathematics, University of Kansas}
\email{purna@math.ku.edu}
\thanks{\emph{Keywords}: deformation of morphisms, surfaces of general type, canonical map, moduli}
\subjclass[2000]{14J29, 14J10,  14B10, 13D10}
\thanks{The first and the second author were partially supported by grants MTM2006--04785 and MTM2009--06964 and by the UCM research group 910772. The first author also thanks the Department of Mathematics of the University of Kansas for its hospitality. The third author thanks the General Research Fund (GRF) of the University of Kansas for partially supporting this research. He also thanks the Algebra Department of the Universidad Complutense de Madrid for its hospitality.}

\begin{abstract}
In this article we study the deformation of finite maps and show how to use this deformation theory to construct
varieties with given invariants in a projective space. Among other things, we prove a criterion that determines when a
finite map can be deformed to a one--to--one map. We use this criterion to construct new
surfaces of general type with birational canonical map, for different $c_1^2$ and $\chi$ (the canonical map of the surfaces we construct is in fact a finite, birational morphism).
Our general results enable us to describe some new components of the moduli of surfaces of general type.
We also find infinitely many moduli spaces $\mathcal M_{(x',0,y)}$ having one component whose general point corresponds to a canonically embedded surface and another component whose general point corresponds to a surface whose canonical map is a degree $2$ morphism.
\end{abstract}

\maketitle

\section*{Introduction}

In this article we address two themes. Firstly, we study the theory of deformations of a morphism
to projective space that is
finite onto its image. This yields a general criterion
that tells us when such a morphism can be deformed
to a degree $1$ morphism (see Theorem~\ref{onetoone2}). Knowing when a finite morphism can be deformed to a morphism of degree $1$ or, even better, to an embedding, is interesting because of its applications
in various contexts.
For instance, it plays a crucial role in smoothing multiple structures on algebraic varieties (see~\cite{Fong},~\cite{GPcarpets},~\cite{Gon},~\cite{GGPropes},~\cite{GGPcarpets} and~\cite{infinitesimal}).
More interestingly, as the method and techniques developed in this article will show, deforming finite morphisms to embeddings or to finite and birational morphisms can be used to prove the existence of varieties of given invariants equipped with a morphism to projective space that is finite and birational into its image.
This leads to the second theme of our article, which illustrates how the deformation techniques we develop give a new method of constructing
surfaces with birational canonical morphism. Mapping the geography of surfaces of general type, i.e., finding surfaces 
of general type with given invariants, is a problem that has been extensively studied by many algebraic geometers; 
as of now, almost all admissible pairs $(\chi, c_1^2)$ are known to be realized by some surface of general type. 
A subtler version of this problem is to know if, for a given pair  $(\chi, c_1^2)$, there exists a surface of general type 
with these invariants, canonically embedded in $\mathbf P^{p_g-1}$, or, at least, if, for a given pair  $(\chi, c_1^2)$, 
there exists a surface with birational canonical map and these given invariants $(\chi, c_1^2)$. The problem of constructing surfaces with birational canonical map  was posed by Enriques (see~\cite{enriques}, chapter VIII, page 284), who called these surfaces \emph{simple canonical surfaces}. In recent years there have been several works dealing with the construction of this kind of surfaces, like the articles of Ashikaga (see~\cite{Ashikaga}) and Catanese (see~\cite{Cat.Hirzebruch.volume},~\cite{Babagge}). Regarding canonically embedded surfaces, one easy source of examples are complete intersections ($(2,4)$ in  $\mathbf P^4$, $(2,2,3)$ in $\mathbf P^5$, etc.). Another easy
way of constructing canonically embedded surfaces (it essentially comes down to using adjunction) is to look at smooth divisors in threefolds such as rational normal scrolls.
An example of this are the surfaces appearing in~\cite[4.5]{AK}; these surfaces produce infinitely many different pairs $(p_g, c_1^2)$ on the Castelnuovo line $c_1^2=3p_g-7$. Apart from this, constructing surfaces with birational canonical map, not to mention
canonical surfaces, is not in general as easy. Our method has the added advantage of producing surfaces whose canonical map is not only birational but also a finite morphism; such a morphism comes closer to being an embedding.

\medskip

The method 
presented here is based upon the following scheme:
\begin{enumerate}
\item[a)] Start with a smooth surface $Y$ having a reasonably simple structure
(typically with geometric genus $0$; see~\cite[Theorem 3.4]{Be}), embedded in a projective space of the desired dimension;
\item[b)] construct a finite canonical cover $\varphi$ of $Y$ having the desired invariants;
\item[c)] deform $\varphi$ to a degree $1$ morphism.
\end{enumerate}

\medskip

Implementing this method successfully rests on having criteria that tell us when the morphism $\varphi$ in b) above can 
be deformed to a degree $1$ morphism. 
The techniques needed to obtain the criteria are developed in the first two sections. In Theorem~\ref{onetoone2}
of Section~\ref{intrinsic.section}  we obtain a criterion that holds in a very general setting. Rather surprisingly 
Theorem~\ref{onetoone2} reduces the question of whether $\varphi$ deforms to a degree $1$ morphism to looking at the 
infinitesimal deformations of $\varphi$ and, eventually, to a cohomological criterion, that needs to be verified only on $Y$. 
Theorem~\ref{onetoone2} is strengthened in~\cite{infinitesimal}, where we give a criterion to say when a finite morphism can be deformed to an embedding. Therefore the criterion in~\cite{infinitesimal} can be used in constructing new embedded varieties with given invariants.

\medskip

Low degree canonical maps are of deep interest, and have a ubiquitous presence in the geometry of algebraic surfaces.
Degree two canonical covers can be viewed as a higher dimensional analogue of hyperelliptic curves.  Of course the situation for surfaces is far more complicated due to the complexity of its  moduli and
the existence of higher degree canonical maps, some of which admit unbounded families in terms of the geometric genus and irregularity (see~\cite{Hor2p_g-4} and~\cite{quad1}).
In Section~\ref{deformations.canonical.general.section} we specialize the theory developed in Section~\ref{intrinsic.section} to the important case of the canonical map $\varphi$ of varieties of general type,  focusing particularly on those surfaces of general type whose canonical map $\varphi$ is a finite morphism of degree $2$.
The general deformation of $\varphi$ exhibits two possible behaviors: either the general deformation is again a finite morphism  of degree $2$ (and in this case, all deformations of $\varphi$ are alike), or the general deformation is a  finite canonical morphism of degree $1$, which is possibly an embedding.
In Theorems~\ref{Psi2=0} and~\ref{Psi2<>0} we give criteria
that say which of the two behaviors occurs in a given case.
These criteria have the following bearing on the moduli of surfaces of general type: they would tell us if an irreducible component contains a ``hyperelliptic'' locus as a proper subset or, on the contrary, if an entire irreducible component is ``hyperelliptic''.  Here by ``hyperelliptic'' locus  we mean the locus parameterizing surfaces of general type whose canonical map is a morphism of degree $2$.

\medskip

In Section~\ref{construct.from.blownup.section} we achieve two goals which are intertwined. The first goal is to illustrate how our method for constructing surfaces with finite, birational canonical morphism works.
 To do this we consider canonical double covers of surfaces $Y$,  where
 $Y$ is $\mathbf P^2$ blown--up at $s$ points in general position, embedded by the system of plane curves of degree $d$ through the $s$ base points. We then apply Theorem~\ref{Psi2<>0} to these covers.
The second goal consists in studying the behavior of the general deformation of the canonical morphism of all surfaces of general type that are canonical double covers of such surfaces $Y$.
To do so we need to use not only Theorem~\ref{Psi2<>0} but also
Theorem~\ref{Psi2=0}.

\medskip

Regarding the first goal mentioned in the previous paragraph, we give now a few more details on how we run our method 
in Section~\ref{construct.from.blownup.section}. Part a) of the method is taken care in Proposition~\ref{blowup.embeddings}. 
Once $Y$ is embedded we go on with Part b). For this we need to construct canonical double covers of $Y$. 
This depends on the existence of smooth divisors on $Y$ in a suitable linear equivalence class, which is $|\omega_Y^{-2}(2)|$.

To determine the existence of a smooth member gets especially harder if the dimension of $|\omega_Y^{-2}(2)|$ is small, for then general arguments do not work.  In fact, understanding the dimension and characteristics of linear systems of plane curves is far from easy (see for instance~\cite{dAH},~\cite{Cop} or~\cite{CilibertoMiranda})
and there are many unsettled questions in this regard. In Lemma~\ref{Harbourne} and  Proposition~\ref{exist.emb.cover}, we give a criterion for $|\omega_Y^{-2}(2)|$ to be base--point--free. As a consequence, we obtain  Corollary~\ref{exist.emb.cover.cor} which settles the existence or non--existence of smooth canonical double covers for most values of $d$ and $s$.

\medskip

Once suitable canonical double covers as above are constructed we go on with Step c) in our method,
which is to determine if
these covers can be deformed to morphisms of degree $1$. Theorem~\ref{Psi2<>0} is a crucial tool in this regard.
As noted before, a consequence  of Theorem~\ref{Psi2<>0}
can be expressed geometrically as the existence of smoothings of certain ``canonical'' non reduced schemes supported on $Y$ (see Corollary~\ref{exist.abstr.carpets} and Theorem~\ref{exist.embedd.carpets}).
To apply Theorem~\ref{Psi2<>0} we need to verify its two crucial conditions 3) and 4), which boil down
to verification of the vanishing and non--vanishing of certain cohomology groups of vector bundles on $Y$. Interestingly, we reduce this to the study of maps of multiplication of global sections of line bundles on  $\mathbf P^2$(see Lemma~\ref{surjectivity} and Proposition~\ref{vanish.normal}).

\medskip

The consequence of all this is the main result of Section~\ref{construct.from.blownup.section}, namely
Theorem~\ref{theorem.construct.surfaces}, which proves the existence of regular surfaces with finite birational canonical morphism and these invariants (see Proposition~\ref{cancover.blownP2.invariants}):

\medskip

\centerline{\vbox{\tabskip=0pt \offinterlineskip
%\tabskip= .25 truecm
\def\tablerule{\noalign{\hrule}}
\halign to 4truecm
%{\valign to125pt}
{\strut
#& \vrule%#
\tabskip=0em plus 3em
%\tabskip= .25 truecm
 %\hskip .2cm
\hfil
#
\hfil  %\hskip .1cm
&
%\hskip .2cm
\hfil#\hfil %\hskip .1cm
& \vrule#&
%\hskip %.2cm
\hfil#\hfil %\hskip %.1cm
& \vrule#&
%\hskip %.2cm
\hfil#\hfil %\hskip %.1cm
& \vrule#&
%\hskip %.2cm
\hfil#\hfil %\hskip %.1cm
& \vrule#\tabskip=0pt
\cr\tablerule
%&&Type
&&%\omit\hidewidth
$p_g$ %\hidewidth
&&%\omit\hidewidth
$q$ %\hidewidth
&&%\omit\hidewidth
 $\chi$  %\hidewidth
&&%\omit\hidewidth
  $c_1^2$  %hidewidth
&\cr\tablerule
%\cr
\tablerule
\tablerule
&&$5$
&&$0$
&&$6$
&&$8$
&\cr
%&\vskip -2truecm&\vskip -2truecm&\vskip -2truecm&\vskip -2truecm&\vskip -2truecm&\vskip -2truecm&\vskip -2truecm&\vskip -2truecm&\vskip -2truecm&\vskip -2truecm&\vskip -2truecm&\vskip -2truecm&\vskip -2truecm&\vskip -2truecm&\cr
\tablerule
%&&2
&&$4$
&&$0$
&&$5$
&&$6$
&\cr
\tablerule
%&&3
&&$7$
&&$0$
&&$8$
&&$16$
&\cr
\tablerule
%&&4
&&$6$
&&$0$
&&$7$
&&$14$
&\cr
\tablerule
%&&5.1
&&$5$
&&$0$
&&$6$
&&$12$
&\cr
\tablerule
&&$8$
&&$0$
&&$9$
&&$24$
&\cr
\tablerule
%&&6.1
&&$7$
&&$0$
&&$8$
&&$22$
&\cr
\tablerule
\noalign{\smallskip}
}}}

\smallskip

\noindent
Except for the case $(p_g,c_1^2)=(5,8)$ and possibly $(p_g,c_1^2)=(7,16)$, surfaces with the above invariants are not complete intersections in projective space. In addition, except for the case $(p_g,c_1^2)=(5,8)$, surfaces with invariants as in the above table cannot be smooth divisors  on $3$--dimensional rational normal scrolls  either (see Proposition~\ref{noin3folds}). This means that surfaces with birational canonical map and invariants $(p_g,c_1^2)=(4,6), (5,12), (6,14), (7,22)$ and $(8,24)$ cannot be constructed by simple methods.

\smallskip

The surfaces with invariants $(p_g,c_1^2)=(4,6)$ were known to Max Noether and Enriques and, together with surfaces with invariants $(p_g,c_1^2)=(5,8)$, were extensively studied by Horikawa (see~\cite{Horikawa.small3} and~\cite{Horikawa.small4}).
In \cite[Theorem 3.2 (a)]{Ashikaga} Ashikaga obtained, by a completely different method, regular surfaces with 
invariants $(p_g,c_1^2)=(7,16), (6,14), (5,12)$. His method only allowed him
to prove the canonical map to be birational but not necessarily a finite morphism. In contrast, our method proves the existence of regular surfaces with invariants  $(p_g,c_1^2)=(7,16), (6,14), (5,12)$  and whose canonical map is a finite birational morphism. Especially interesting is the case $(p_g,c_1^2)=(5,12)$, for the existence of smooth surfaces in $\mathbf P^4$ is severely constrained. Specifically, this constraint is given by the double point formula, which implies that the only smooth, canonically embedded surfaces in $\mathbf P^4$ are the ones with $(p_g,c_1^2)=(5,8)$, which are in fact complete intersections of type $(2,4)$. In this context we show the existence of surfaces with $(p_g,c_1^2)=(5,12)$ whose canonical map, which cannot be an embedding because of the double point formula, comes as close to being an embedding as possible.  For the other invariants in the above table (that is,
$(p_g,c_1^2)=(8,24), (7,22)$) surfaces with birational canonical map were not previously known to exist, to the best of our understanding.  The surfaces we construct have finite (birational) canonical morphism and, rather surprisingly,
are regular. In fact Debarre proved (see~\cite{Debarre})
that  a minimal irregular surface of general type should satisfy the inequality $c_1^2 \geq 2p_g$ so one might expect  surfaces above the line $c_1^2 = 2p_g$ to be more likely irregular.
However, the surfaces we construct with invariants $(p_g,c_1^2)=(7,16), (6,14), (5,12), (8,24), (7,22)$ are
regular despite satisfying $c_1^2 >  2p_g$. We also point out that these surfaces are not only regular, but also simply connected (see Remark~\ref{simply.connected}.)

\medskip
The second goal of Section~\ref{construct.from.blownup.section} is to determine the behavior of the general deformation of the canonical morphism of canonical double covers of surfaces $Y$ which are $\mathbf P^2$ blown--up at $s$ points in general position and embedded by the system of plane curves of degree $d$ through the $s$ base points.
This is done in Theorems~\ref{theorem.construct.surfaces} and~\ref{2to1deformsto2to1}, which say that,  in all but finitely many cases, the canonical double cover of $Y$ cannot be deformed to a morphism of degree $1$. The finitely many exceptions are precisely those cases that allow us to construct the surfaces with finite and birational canonical morphism mentioned in the previous paragraphs. Thus most surfaces $X$ obtained as canonical double covers of surfaces $Y$
live in a ``hyperelliptic''
component of the moduli, in the sense explained earlier.  Continuing in the same strain,
 we show in~\cite{Hirzebruch} that a
 canonical double cover of any embedded minimal rational surface always deforms to a canonical double cover.
 Thus Theorem~\ref{2to1deformsto2to1} and the results of~\cite{Hirzebruch}  demonstrate the existence of ``hyperelliptic'' 
 components for infinitely many moduli spaces.
 The above does not mean that this is the generic situation for the ``hyperelliptic'' locus of moduli spaces of surfaces of  general type.
 For instance, the construction of~\cite[4.5]{AK} and its generalization Lemma~\ref{divisor.scroll} produce components with 
 a proper ``hyperelliptic locus'' for infinitely many moduli spaces. The existence of ``hyperelliptic'' 
 components`` and of components with 
 a proper ``hyperelliptic locus'', as shown by Theorems~\ref{theorem.construct.surfaces} and~\ref{2to1deformsto2to1}, the results in~\cite{Hirzebruch}, the construction of~\cite[4.5]{AK} and Lemma~\ref{divisor.scroll}, bears witness to the complexity of 
 the moduli of surfaces of general type.

\medskip
Finally, in Section~\ref{moduli.section} we compute the invariants
and we study the moduli components of the surfaces of general type $X$ constructed in Section~\ref{construct.from.blownup.section}. In this regard, besides showing the curious phenomena described in the previous paragraph, we compute the dimension $\mu$ of the moduli components of the surfaces $X$.
When $X$ is as in Theorem~\ref{theorem.construct.surfaces}, the moduli component of $X$ has a proper ``hyperelliptic'' locus, whose dimension $\mu_2$ we also compute in Proposition~\ref{prop.moduli.2to1deformsto1to1}, obtaining the following table:

\medskip

 \centerline{\vbox{\tabskip=0pt \offinterlineskip
%\tabskip= .25 truecm
\def\tablerule{\noalign{\hrule}}
\halign to 3.65truecm
%{\valign to125pt}
{\strut
#& \vrule#
\tabskip=0em plus 3em
%\tabskip= .25 truecm
& \hskip .25cm
\hfil #
\hfil  \hskip .05cm
& \vrule #
& \hskip .3cm
\hfil #
\hfil  \hskip .05cm
& \vrule \vrule \vrule #
&\hskip .3cm
\hfil# \hskip .05cm
&  \vrule #
& \hskip .25cm
\hfil #
\hfil  \hskip .05cm
& \vrule #&
\hskip .15cm \hfil# %\hskip .05truecm
\hfil & \vrule#\tabskip=0pt
\cr\tablerule
%&&Type
&&%\omit\hidewidth
$d$ %\hidewidth
&& %\omit\hidewidth
$s$ %\hidewidth
&& %\omit\hidewidth
$\mu$ %\hidewidth
&&%\omit\hidewidth
$\mu_2$ %\hidewidth
&\cr\tablerule
%\cr
\tablerule
\tablerule
&&$3$
&&$5$
&&$44$
&&$42$
&\cr
%&\vskip -2truecm&\vskip -2truecm&\vskip -2truecm&\vskip -2truecm&\vskip -2truecm&\vskip -2truecm&\vskip -2truecm&\vskip -2truecm&\vskip -2truecm&\vskip -2truecm&\vskip -2truecm&\vskip -2truecm&\vskip -2truecm&\vskip -2truecm&\cr
\tablerule
%&&2
&&$3$
&&$6$
&&$38$
&&$34$
&\cr
\tablerule
%&&3
&&$4$ %\ *$ \hskip -.2truecm
&&$8$
&&$48$
&&$47$
&\cr
\tablerule
%&&4
&&$4$ % \ *$ \hskip -.2truecm
&&$9$
&&$42$
&&$39$
&\cr
\tablerule
%&&5.1
&&$4$ %\hskip -.15truecm
&&$10$
&&$36$
&&$31$
&\cr
\tablerule
&&$5$ %\hskip -.15truecm
&&$13$
&&$42$
&&$40$
&\cr
\tablerule
%&&6.1
&&$5$ %\hskip -.15truecm
&&$14$
&&$36$
&&$32$
&\cr
\tablerule
\noalign{\smallskip}}}}

\medskip

\noindent In particular we recover results of Horikawa (see~\cite{Horikawa.small3} and~\cite{Horikawa.small4}) about the moduli components corresponding to surfaces $X$ with invariants $(p_g, c_1^2)=(4,6)$ and $(p_g,c_1^2)=(5,8)$  that appear in  Theorem~\ref{theorem.construct.surfaces}.

\medskip

In the last part of Section~\ref{moduli.section} we combine Theorem~\ref{2to1deformsto2to1} and Lemma~\ref{divisor.scroll} 
to obtain Theorem~\ref{2components},  where we show,
for infinitely many values $(p_g,0,c_1^2)$, that the moduli space of surfaces of general type
has one component which is ``totally hyperelliptic'' and another component whose general point corresponds to a canonically embedded surface. In many cases this latter component has a proper ``hyperelliptic locus'', see Lemma~\ref{divisor.scroll}.
This bears further evidence to the complexity of the moduli of surfaces of general type, compared for instance to the moduli of curves.
It would be nice to understand the more philosophical, deeper reasons behind these differences in behavior for the components of the moduli of surfaces.

\section{Sufficient conditions to deform finite morphisms to morphisms of degree $1$}\label{intrinsic.section}

The purpose of this section is to give criteria to assure when a morphism from a smooth variety to projective space, finite onto a smooth image, can be deformed to a degree $1$ morphism.

\medskip

\begin{noname}\label{setup}
{\bf  Notation and conventions:} {\rm Throughout this article, unless otherwise stated,
we will use the following notation and conventions:
\begin{enumerate}
\item  We will work over an algebraically closed field $\mathbf k$ of characteristic $0$.
\item  $X$ and $Y$ will denote smooth, irreducible projective varieties.
\item  $i$ will denote a projective embedding $i: Y \hookrightarrow \mathbf P^N$. In this case, $\mathcal I$  will denote the ideal sheaf of $i(Y)$ in $\mathbf P^N$. Likewise, we will often abridge $i^*\mathcal O_{\mathbf P^N}(1)$ as $\mathcal O_Y(1)$.
\item $\pi$ will denote a finite morphism $\pi: X \longrightarrow Y$ of degree $n \geq 2$; in this case, $\mathcal E$ will denote  the trace--zero module of $\pi$ ($\mathcal E$ is a vector bundle on $Y$ of rank $n-1$).
\item $\varphi$ will denote a projective morphism $\varphi: X \longrightarrow \mathbf P^N$ such that $\varphi= i \circ \pi$.
\end{enumerate}}
\end{noname}

\noindent We introduce a homomorphism defined in~\cite[Proposition 3.7]{Gon}:

\begin{proposition}\label{morphism.miguel}
There exists a homomorphism
\begin{equation*}
 H^0(\mathcal N_\varphi) \overset{\Psi}\longrightarrow \mathrm{Hom}(\pi^*(\mathcal I/\mathcal I^2), \mathcal O_X),
\end{equation*}
that appears when taking cohomology on the commutative diagram~\cite[(3.3.2)]{Gon}. Since
\begin{equation*}
\mathrm{Hom}(\pi^*(\mathcal I/\mathcal I^2), \mathcal O_X)=\mathrm{Hom}(\mathcal I/\mathcal I^2, \pi_*\mathcal O_X)=\mathrm{Hom}(\mathcal I/\mathcal I^2, \mathcal O_Y) \oplus \mathrm{Hom}(\mathcal I/\mathcal I^2, \mathcal E)
\end{equation*}
the homomorphism $\Psi$ has two components
\begin{eqnarray*}
H^0(\mathcal N_\varphi) & \overset{\Psi_1}  \longrightarrow & \mathrm{Hom}(\mathcal I/\mathcal I^2, \mathcal O_Y) \cr
H^0(\mathcal N_\varphi) & \overset{\Psi_2}  \longrightarrow & \mathrm{Hom}(\mathcal I/\mathcal I^2, \mathcal E).
\end{eqnarray*}
\end{proposition}

\begin{proposition}\label{onetoone}
Let $T$ be
a smooth irreducible algebraic curve
with a distinguished point $0$.
Let $\Phi: \mathcal X \longrightarrow \mathbf P^N_T$ be a flat family of morphisms over $T$ (i.e., $\Phi$ is a $T$--morphism for which $\mathcal X \longrightarrow T$ is proper, flat and surjective) such that
\begin{enumerate}
\item $\mathcal X$ is irreducible and reduced;
\item $\mathcal X_t$ is smooth, irreducible and projective for all $t \in T$;
\item $\mathcal X_0=X$ and $\Phi_0=\varphi$.
\end{enumerate} Let $\Delta$ be the first infinitesimal neighborhood of $0$ in $T$. Let $\tilde X$ and $\tilde \varphi$ be the pullbacks of $\mathcal X$ and $\Phi$ over $\Delta$ and let $\nu$ be the element of $H^0(\mathcal N_\varphi)$ corresponding to $\tilde \varphi$.
If  $\Psi_2(\nu) \in \ \mathrm{Hom}(\mathcal I/\mathcal I^2, \mathcal E)$ is a homomorphism of rank $k > n/2-1$, then, shrinking $T$ if necessary, $\Phi_t$ is finite and one--to--one for any $t \in T, t \neq 0$.
\end{proposition}

\begin{proof}
Recall (see~\cite[Proposition 2.1.(1)]{Gon}) that $\Psi_2(\nu)$ corresponds to a pair $(\tilde Y, \tilde i)$ where $\tilde Y$ is a rope on $Y$ with conormal bundle $\mathcal E$ (therefore, in particular, $\tilde Y$ is a rope of multiplicity $n$ on $Y$) and $\tilde i: \tilde Y \longrightarrow \mathbf P^N$ extends $i$.
Then $\mathcal E$ is both the conormal bundle of $\tilde Y$ and the trace--zero module of $\pi$.
Let $L=\varphi^*\mathcal O_{\mathbf P^N}(1)$ and $\mathcal L=\Phi^*\mathcal O_{\mathbf P^N_T}(1)$.  Let $\mathcal Y=\Phi(\mathcal X)$ and define $\Pi: \mathcal X \longrightarrow \mathcal Y$ so that $\Phi$ factors through $\Pi$.
%
%\noindent
Shrinking $T$ if necessary we may assume that $\mathcal Y_t$ is the image of $\Phi_t$ for all $t \in T, t \neq 0$. Then
\begin{equation}\label{inequalities}
n \, \mathrm{deg} \, Y = \mathrm{deg}\, L=\mathrm{deg}\, \mathcal L_t=\mathrm{deg} \, \Pi_t\cdot \mathrm{deg}\, \mathcal Y_t=\mathrm{deg}\, \Pi_t\cdot \mathrm{deg}\, \mathcal Y_0 \geq \mathrm{deg}\, \tilde i(\tilde Y) \cdot \mathrm{deg}\, \Pi_t =(k+1)\mathrm{deg}\, Y \cdot \mathrm{deg}\, \Pi_t.
\end{equation}
To justify~\eqref{inequalities} we recall that, since $\mathcal X$ is a flat family, $\mathrm{deg}\,L=\mathrm{deg}\,\mathcal L_t$
and $\mathrm{deg}\,\mathcal Y_0 \geq \mathrm{deg}\,\tilde i(\tilde Y)$, since $\tilde i(\tilde Y) \subset \mathcal Y_0$ because $(\mathrm{im}\,\tilde \varphi)_0=\tilde i(\tilde Y)$ by~\cite[Theorem 3.8.(1)]{Gon}.
We also use that $k+1$ is the multiplicity of $\tilde i(\tilde Y)$. Then from~\eqref{inequalities} it follows that
\begin{equation*}
n \geq (k+1)  \mathrm{deg}\,\Pi_t
\end{equation*}
which, together with the hypothesis $k > n/2-1$, implies that $\mathrm{deg}\,\Pi_t < 2$. Finally, $\Pi_t$ is finite because so is $\pi$.
\end{proof}

\noindent We will use Proposition~\ref{onetoone} to obtain Theorem~\ref{onetoone2}, which is a criterion to know when a projective morphism from a smooth variety, finite onto a smooth image, can be deformed to a degree $1$ morphism.

\begin{theorem}\label{onetoone2}
Let $\tilde \varphi: \tilde X \longrightarrow \mathbf P^N_\Delta$ be a first order infinitesimal deformation of $\varphi$ ($\Delta$ is Spec$(\mathbf k[\epsilon]/\epsilon^2)$) and let $\nu$ be the class of $\tilde \varphi$ in $H^0(\mathcal N_\varphi)$. If
\begin{enumerate}
\item the homomorphism $\Psi_2(\nu) \in \ \mathrm{Hom}(\mathcal I/\mathcal I^2,\mathcal E)$ has rank $k > n/2-1$; and
\item  there exists an algebraic formally semiuniversal deformation of $\varphi$ and $\varphi$ is unobstructed,
\end{enumerate}
then there exists a flat family of morphisms, $\Phi: \mathcal X \longrightarrow \mathbf P^N_T$ over $T$, where $T$ is
a smooth irreducible algebraic curve with a distinguished point $0$,
such that
\begin{enumerate}
\item[(a)] $\mathcal X_t$ is a smooth, irreducible, projective variety;
\item[(b)] the restriction of $\Phi$ to the first infinitesimal neighborhood of $0$ is $\tilde \varphi$ (and hence $\Phi_0=\varphi$); and
\item[(c)] for any $t \in T$, $t \neq 0$, $\Phi_t$ is finite and one--to--one onto its image in $\mathbf P^N_t$.
\end{enumerate}
\end{theorem}

\begin{proof}
There exists a smooth algebraic variety
$M$ which is the base of the algebraic formally semiuniversal deformation
of $\varphi$
and $H^0(\mathcal N_\varphi)$ is the tangent space of $M$ at $[\varphi]$.
Then $\nu$ represents a tangent vector to $M$ at $[\varphi]$.  Thus there is a smooth algebraic curve
$T \subset M$ passing through $[\varphi]$ and a family of morphisms over $T$,  $\Phi: \mathcal X \longrightarrow \mathbf P^N_T$ satisfying (1), (2) and (3) in Proposition~\ref{onetoone} and $\Phi_\Delta=\tilde \varphi$. Now, since $\Psi_2(\nu)$ is a  homomorphism of rank $k > n/2-1$, Proposition~\ref{onetoone} implies (c).
\end{proof}

\begin{remark}
 {\rm The criterion obtained in Theorem~\ref{onetoone2} holds in a very general setting.  Indeed, the requirement on $\varphi$ to possess an algebraic formally semiuniversal deformation is satisfied if $\varphi$ is non--degenerate (for the definition of non--degenerate morphisms, see~\cite[p. 376]{Hor1} or~\cite[Definition 3.4.5]{Ser}) and if ${X}$ possesses a formal semiuniversal deformation which is effective, since in such a case the formal semiuniversal deformation of $\varphi$ is also effective and therefore algebraizable by Artin's algebraization theorem (see \cite{Artin.alg}). In particular, if $\varphi$ is the canonical map of a smooth variety of general type $X$ with ample and base--point--free canonical bundle, $\varphi$ possesses an algebraic formally semiuniversal deformation (see Lemma~\ref{effective-1}).}
\end{remark}

\section{Deformations of the canonical morphism of a variety of general type with ample canonical line bundle}\label{deformations.canonical.general.section}

In this section we will apply the theory developed in Section~\ref{intrinsic.section} to the study of canonical morphisms of varieties of general type.
The first lemma gives conditions to tell when a morphism is unobstructed. Recall that the unobstructedness of $\varphi$ was required in Theorem~\ref{onetoone2}. Lemma~\ref{varphi.iff.X.unobstructed} applies nicely to the situation of canonical morphisms that we study later on in this section.

\begin{lemma}\label{varphi.iff.X.unobstructed}
Let $X$ be a smooth regular projective variety.
Let $L$ be a globally generated line bundle on $X$ such that $H^1(L)=0$ and $L$ admits a lifting to any first--order infinitesimal deformation of $X$ (the latter condition is equivalent to the vanishing of the map $H^1(\ST_{X}) \to H^2(\SO_{X})$ induced by the cohomology class $c_1(L) \in H^1(\Omega_{X})$ via cup product and duality). Let $X \xrightarrow{\psi} \PP^N$ be the map induced by the complete linear series of $L$. Assume $\psi$ is non--degenerate. Then the forgetful morphism $\mathrm{Def}_{\psi} \to \mathrm{Def}_{X}$, where $\mathrm{Def}_{\psi}$ is the functor of infinitesimal deformations of $\psi$ with fixed target and $\mathrm{Def}_{X}$ is the functor of infinitesimal deformations of $X$, is smooth. As a consequence, $\psi$ is unobstructed if and only if $X$ is unobstructed.
\end{lemma}

\begin{proof} Since each functor admits a semiuniversal formal element,
to see that the forgetful morphism
\begin{equation}\label{forgetful.morphism}
\mathrm{Def}_{\psi} \to \mathrm{Def}_{X}
\end{equation}
is smooth it is enough, by general deformation theory, to check that the differential of~\eqref{forgetful.morphism} is surjective and that $\mathrm{Def}_{\psi}$ is less obstructed than $\mathrm{Def}_{X}$.

\smallskip

\noindent
To see that the differential of~\eqref{forgetful.morphism} is surjective recall that by hypothesis $L$ lifts to any given first--order deformation of $X$ and $H^1(L)=0$. Then, bases for the global sections of $L$ can be lifted also to bases for the global sections of any lifting $\WL$ of $L$ to any first--order deformation of $X$. Thus $\psi$ can be lifted to any infinitesimal deformation of $X$, so
the differential of~\eqref{forgetful.morphism}
is surjective.

\smallskip

\noindent Now we see that $\mathrm{Def}_{\psi}$ is less obstructed than $\mathrm{Def}_{X}$. For this we will see first that the functor $\mathrm{Def}_{\psi}^{f}$ of infinitesimal deformations of $\psi$ with fixed domain and target is unobstructed (from now on we will abridge the phrase ``with fixed domain and target''  by \emph{wfdat}, as done in~\cite[p. 157]{Ser}). By general deformation theory, the unobstructedness of $\mathrm{Def}_{\psi}^{f}$ is equivalent to checking that, for any small extension of Artin rings
\begin{equation*}
 0 \to (t) \to \WA \to A \to 0,
\end{equation*}
 any wfdat deformation of $\psi$ over $A$  lifts to a wfdat deformation over $\WA$  (see e.g. \cite[pg. 48]{Ser}).
To see this, note that the condition $H^1(\SO_{X})=0$ implies that the local Picard functor of $(X, L)$ is constant (see e.g. \cite[Theorem 3.3.1]{Ser}). 
This means
that the only lifting of $L$ to any trivial infinitesimal deformation of $X$ is the trivial one.
Moreover, it follows from the sequence
\begin{equation*}
 0 \to (t)\otimes L \to \WA \otimes L \to A \otimes L \to 0
\end{equation*}
and from $H^1(L)=0$ that global sections of the trivial deformation of $L$ over $A$ lift to global sections of the trivial deformation of $L$ over $\WA$. This proves that any wfdat deformation of $\psi$ over $A$  lifts to a wfdat deformation over $\WA$ 
so, as explained before, $\mathrm{Def}_{\psi}^{f}$
is unobstructed.

\smallskip

\noindent
Now, note that the differential of~\eqref{forgetful.morphism} is the same as the connecting homomorphism
\begin{equation*}
H^0(\SN_{\psi}) \to H^1(\ST_{X}),
\end{equation*}
so the surjectivity of the differential of~\eqref{forgetful.morphism} gives us the injectivity at the left hand side of the exact sequence
\begin{equation}\label{obstr.space.seq}
0 \to H^1(\psi^{*}\ST_{\PP^N}) \to H^1(\SN_{\psi}) \to H^2(\ST_{X}).
\end{equation}
The spaces appearing in sequence~\eqref{obstr.space.seq} are, from left to right, obstruction spaces for the functors $\mathrm{Def}_{\psi}^{f}$, $\mathrm{Def}_{\psi}$ and $\mathrm{Def}_{X}$. Since $\mathrm{Def}_{\psi}^{f}$ is unobstructed we conclude that the obstructions for $\mathrm{Def}_{\psi}$ actually live in
\begin{equation*}
H^1(\SN_{\psi})/H^1(\psi^{*}\ST_{\PP^N}).
\end{equation*}
Thus it follows from sequence~\eqref{obstr.space.seq} that $\mathrm{Def}_{\psi}$ is less obstructed than $\mathrm{Def}_{X}$.
\end{proof}

\noindent Now we state some results about deformations of varieties of general type and its canonical maps, needed to apply Theorem~\ref{onetoone2} and to prove Theorems~\ref{Psi2=0} and~\ref{Psi2<>0}. The first of this results is the following corollary of Lemma~\ref{varphi.iff.X.unobstructed}:

\begin{cor}\label{canonical.remark} Let $X$ be a smooth projective variety of dimension $m \geq 2$.
\begin{enumerate}
\item If the canonical bundle of $X$ is base--point--free and $H^1(\SO_{X})=H^{m-1}(\SO_{X})=0$, then the canonical morphism of $X$ is unobstructed if and only if $X$ is unobstructed.
\item If $X$ is regular, the canonical bundle of $X$ is ample and the pluricanonical bundles of $X$ are base--point--free, then the pluricanonical morphims of $X$ are unobstructed if and only if $X$ is unobstructed.
\item In particular, if $X$ is a regular, smooth projective surface of general type with ample and base--point--free canonical  bundle, then the canonical morphism or any of the pluricanonical morphisms of $X$ is unobstructed if and only if $X$ is unobstructed.
\end{enumerate}
\end{cor}

\begin{proof}
The result will follow from Lemma~\ref{varphi.iff.X.unobstructed} once we check that the hypotheses of the lemma hold.  For (1), we have by hypothesis that the canonical bundle is base--point--free and that $H^1(\SO_{X})=0$. Also $H^1(\omega_{X})=H^{m-1}(\SO_{X})^{\vee}=0$ and the canonical bundle lifts to the relative dualizing sheaf. For (2) the argument is similar except that we use the Kodaira vanishing to obtain $H^1(\omega_{X}^{\otimes l})=0$. Finally (3) is a straight--forward consequence of (1) and (2).
\end{proof}

\begin{lemma}\label{effective-1}
If  $X$ is a smooth projective variety of general type with ample canonical bundle, then $X$ has an algebraic formally universal deformation. If in addition the canonical bundle of $X$ is base--point--free,
then the canonical or $l$--pluricanonical morphism ${X} \xrightarrow{\psi} \PP^N$ has an algebraic formally universal deformation.
\end{lemma}

\begin{proof}
The first assertion is well known: the existence of an
algebraic formally semiuniversal deformation for $X$ follows from
Grothendieck's existence theorem (see \cite[Theorem 5.4.5]{EGA3-1} or \cite[Theorem 2.5.13]{Ser})
and Artin's algebraization theorem \cite{Artin.alg}. This deformation is formally universal because $H^0(\ST_{X})=0$.

\noindent
On the other hand, if the canonical map of $X$ is a morphism, then general deformation theory and  Artin's algebraization theorem
imply that the canonical or $l$--pluricanonical morphism ${X} \xrightarrow{\psi} \PP^N$ has an algebraic formally semiuniversal space of deformations which is, in fact, formally universal because $H^0(\ST_{X})=0$ implies not only that $X$ has no infinitesimal automorphisms, but also that $\psi$ has no infinitesimal automorphisms.
\end{proof}

\begin{lemma}\label{canonical.deforms.canonical}
Let $X$ be a smooth
regular projective variety of general type with ample and base--point--free canonical divisor. Let $\psi: X \longrightarrow \mathbf P^N$ be the canonical morphism of $X$. Let $\Phi: \mathcal X \longrightarrow \mathbf P_T^N$ be a deformation of $\psi$ over a smooth curve $T$. Then, after maybe shrinking $T$, $\Phi_t: \mathcal X_t \longrightarrow \mathbf P_t^N$ is the canonical morphism of $\mathcal X_t$ for all $t \in T$.
\end{lemma}

\begin{proof}
We first prove that $\Phi^*\mathcal O_{\mathbf P_t^N}(1)=\omega_{\mathcal X_t}$. Let $L_t=\Phi^*\mathcal O_{\mathbf P_t^N}(1)$ and let $K_t$ be the canonical divisor of $\mathcal X_t$.
Since $R^2 p_*\mathbf Z$ is locally constant (where $p:\mathcal X \to T$),
if $K_0-L_0=0$, then $c_1(K_0-L_0)=0$, and so, $c_1(K_t-L_t)=0$. Then $h^1(\mathcal O_{X})=0$
yields $K_t-L_t=0$.

\noindent
Now recall that $p_g$ is a deformation invariant, so
$h^0(\omega_{X})=h^0(\omega_{\mathcal X_t})$ for all $t \in T$, so $\Phi_t$ is induced by the complete linear series of $\omega_{\mathcal X_t}$.
\end{proof}

\noindent
Finally, for the results of the remaining of the article, we will need to understand how the normal bundle of a double cover is. So we do in the next lemma, which is stated in slightly greater generality, since it considers simple cyclic covers of arbitrary degree.

\begin{lemma}\label{normal.pi} Let $X$, $Y$ and $\pi$ be as in~\ref{setup}.
Assume in addition that $\pi$ is a simple cyclic cover branched along a smooth divisor $B$ of $Y$ and determined by a line bundle $\SL$ on $Y$ (see~\cite[I.17]{BHPV}). Let $R \subset X$ be the ramification divisor of $\pi$  (see~\cite[I.16]{BHPV} for the definition of ramification divisor).
Then the normal bundle $\SN_{\pi}$ of $\pi$ satisfies
\begin{equation}\label{normalpi}
\SN_{\pi}=\SO_R(R)\otimes \pi^*\SL= \pi^*\SL^{n} \otimes \SO_R=\SO_R(\pi^* B).
\end{equation}
\end{lemma}

\begin{proof}
Since $X$ and $Y$ are smooth, from the sequence defining $\SN_{\pi}$,
we see that
\begin{equation}\label{ext.uno}
\SN_{\pi}= \mathcal{E}xt^1(\Omega_{X/Y}, \SO_X).
\end{equation}
Since $\pi$ is a simple cyclic cover (i.e., an $n$-cyclic covering in the notation of~\cite{BHPV}), Lemma I.16.1 of~\cite{BHPV} implies
\begin{equation}\label{normal.r}
\SO_R(R)= \pi^*\SL^{\otimes (n-1)} \otimes \SO_R.
\end{equation}
%\noindent
It is also easy to check that for a simple cyclic cover
\begin{equation}\label{claim}
\Omega_{X/Y}= \pi^*\SL^{-1} \otimes \SO_R.
\end{equation}
Then, from \eqref{ext.uno}, \eqref{normal.r}, \eqref{claim} and the formula $\SO_R(R)=  \mathcal{E}xt^1(\SO_R, \SO_X)$  we obtain \eqref{normalpi}.
\end{proof}

\begin{theorem}\label{Psi2=0}
With the notation of~\ref{setup} and of Proposition~\ref{morphism.miguel}, let $X$ be a smooth variety of general type of dimension $m \geq 2$ with ample and base--point--free canonical bundle, let $\varphi$ be its canonical morphism and assume that the degree of $\pi$ is $n=2$.
Assume furthermore that
\begin{enumerate}
\item $h^1(\mathcal O_Y)=h^{m-1}(\mathcal O_Y)=0$ (in particular, $Y$ is regular);
\item $h^1(\mathcal O_Y(1))=h^{m-1}(\mathcal O_Y(1))=0$;  
\item $h^0(\omega_Y(-1))=0$;
\item $h^1(\omega_Y^{-2}(2))=0$;
\item the variety $Y$ is unobstructed in $\mathbf P^N$ (i.e., the base of the universal deformation space of $Y$ in $\mathbf P^N$ is smooth); and
\item $\Psi_2 = 0$.
\end{enumerate}
Then
\begin{enumerate}
\item[(a)] $X$ and $\varphi$ are unobstructed, and
\item[(b)] any deformation of $\varphi$ is a (finite) canonical morphism of degree $2$. Thus the canonical map of a variety  corresponding to a general point of the component of $X$ in its moduli space is a finite morphism of degree $2$.
\end{enumerate}
\end{theorem}

\begin{proof}
Let $p: \SY \hookrightarrow \PP_{U}^N \to U$ be an algebraic formally universal embedded deformation of $Y$.
We can, by assumption, take both $U$ and the total family $\SY$ smooth.

\smallskip
\noindent
Let us denote $\mathbb L=\omega_{\SY / U}^{-1} \otimes \SO_{\SY}(1)$, where $\omega_{\SY / U}$ is the relative dualizing sheaf on $\SY$.

\smallskip
\noindent
Relative duality implies
\begin{equation}\label{relative.duality}
\pi_* \omega_X = (\pi_* \SO_X)^\vee \otimes \omega_Y= \omega_Y \oplus (\SE^{-1} \otimes \omega_Y).
\end{equation}
Since $\varphi$ is the canonical morphism of $X$ and factors through $\pi$, we have also
\begin{equation}\label{pushdown.canonical}
\pi_* \omega_X = \mathcal O_Y(1) \oplus  \SE(1),
\end{equation}
so taking the determinant on~\eqref{relative.duality} and~\eqref{pushdown.canonical} we conclude
\begin{equation}\label{branch.of.canonical.double.cover}
 \SE^{-2}=\omega_Y^{-2}(2).
\end{equation}
In fact,  \eqref{relative.duality}, \eqref{pushdown.canonical}, \eqref{branch.of.canonical.double.cover} and the connectedness of $X$ imply
that the hypothesis $h^0(\omega_Y(-1))=0$ is equivalent to
\begin{equation*}\label{trace.omega.minus.one}
\SE = \omega_Y(-1).
\end{equation*}
From \eqref{branch.of.canonical.double.cover} we also see that the branch locus $B$ of $\pi$ is a divisor 
in $|\omega_Y^{-2}(2)|$ and, since $X$ is smooth,  $B$ is in fact a smooth
divisor in $|\omega_Y^{-2}(2)|$.

\noindent
By the hypothesis $h^1(\omega_Y^{-2}(2))=0$, we can
assume $h^1((\mathbb L^{\otimes 2})_u)=0$ for any $u \in U$ and $p_*(\mathbb L^{\otimes 2})$ is a free sheaf on $U$ of rank $M+1=h^0(\SO_Y(B))$.

\smallskip
\noindent
Let $\PP(p_*(\mathbb L^{\otimes 2})) \to U$ be the projective bundle associated to $p_*(\mathbb L^{\otimes 2})$. Fix a basis $p_*(\mathbb L^{\otimes 2}) = \SO_U s_0 \oplus \cdots \oplus \SO_U s_M$. Let $X_0, \ldots, X_M \in H^0(\SO_{\PP(p_*(\mathbb L^{\otimes 2}))} (1))=H^0((p_*(\mathbb L^{\otimes 2}))^{\vee})$ be its dual basis. On $\SY \times_U \PP(p_*(\mathbb L^{\otimes 2}))$, we consider the divisor $\SB$ defined by the section $X_0 s_0 + \cdots + X_M s_M$.

\noindent
Let $\SY \times_U \PP(p_*(\mathbb L^{\otimes 2})) \xrightarrow{q} \SY$ denote the projection. Let $t$ denote the tautological section on the total space of $q^*\mathbb L$. Then we can construct a relative double cover
\begin{equation*}\label{universal.double.cover}
\xymatrix@C-25pt{
(t^2- (X_0s_0 + \cdots + X_M s_M))_{0}\ar@{=}[r] & \SX \ar[d] \ar@{^{(}->}[rrrr]& && &
 q^* \mathbb L \ar[lllld]\\
& \SY \times_U \PP(p_*(\mathbb L^{\otimes 2}))&&&&,
}
\end{equation*}
whose fiber over any point $(u, [r]) \in \PP(p_*(\mathbb L^{\otimes 2})) $, with $u \in U$ and $r \in H^0((\mathbb L^{\otimes 2})_u)$, is the double cover $\SX_{(u,[r])} \to \SY_u$ defined by the divisor $\SB_{\mid \SY_u}=(r)_0 \in | \omega_{\SY_u}^{-2} \otimes \SO_{\SY_u}(2) |$. In fact, we restrict the construction to the open set $V \subset \PP(p_*(\mathbb L^{\otimes 2}))$, where the divisors $\SB_{\mid \SY_u}$ are smooth,
 in order to obtain integral, smooth, double covers $\SX_{(u,[r])} \to \SY_u$. The open set $V$ contains the point $(u_0, [r_0])$ that corresponds to the pair $(Y,B)$, so we can assume
$V$ maps surjectively onto $U$.

\smallskip
\noindent
Let $\Phi$ denote the composite map
\begin{equation*}
\SX \to \SY \times_{U} V \hookrightarrow \PP_{V}^N.
\end{equation*}
Then $\Phi$ is an algebraic deformation of $X \xrightarrow{\varphi} \PP^N$. We denote with the same symbols the associated formal family at the point $(u_0, [r_0]) \in V$.

\noindent
Let $\mathrm{Def}_{(Y,B)}$ be the functor of deformations of the pair $(Y \subset \PP^N, B \in | \omega_Y^{-2}(2) |)$. The construction above defines a map
\begin{equation}\label{map.F}
\mathrm{Def}_{(Y,B)} \xrightarrow{F} \mathrm{Def}_{\varphi}.
\end{equation}
The family $(\SY \times_U V \to V, \SB, (u_0, [r_0]) \in V)$ is an algebraic formally universal deformation of $(Y \subset \PP^N, B \in | \omega_Y^{-2}(2) |)$. So the functor $\mathrm{Def}_{(Y,B)}$ is prorepresentable and smooth.

\smallskip
\noindent
Assume we prove that the differential $\mathrm{d}F$ is an isomorphism. Then $F$ and $\mathrm{Def}_{\varphi}$ would be smooth. Therefore the family $(\SX \xrightarrow{\Phi} \PP_V^N \to V, (u_0, [r_0])\in V)$ is an algebraic formally semiuniversal deformation for $\varphi$ and $\mathrm{Def}_{\varphi}$ is unobstructed. Then the fact that $\SX \xrightarrow{\Phi} \PP_V^N$ is semiuniversal for $\varphi$ and  Lemma~\ref{canonical.deforms.canonical} mean that $\varphi$ deforms always to a (canonical) $2:1$ morphism.

\smallskip
\noindent
By hypothesis, $\omega_X$ is
base--point--free and, by conditions (1) and (2),
$H^1(\mathcal O_X)=H^{m-1}(\mathcal O_X)=0$. Then, since
$\varphi$ is unobstructed it follows from Corollary~\ref{canonical.remark}, (1),
that $X$ is also unobstructed and therefore the point of $X$ in its moduli space belongs to a unique component of the moduli.
Then the general point of this component has a canonical map which is a finite morphism of degree $2$.

\medskip
\noindent
To complete the proof it remains to see that $\mathrm{d}F$ is an isomorphism. Here we will use the hypothesis $\Psi_2 = 0$.

\noindent
We will compute the tangent space to $V$ at $(u_0,[r_0])$ by considering the fiber at $(u_0,[r_0])$ of the sequence
\begin{equation*}
 0 \to \ST_{V/U} \to \ST_V \to \nu^* \ST_U \to 0,
\end{equation*}
associated to the projection $V \xrightarrow{\nu} U$.
This way, since $H^1(\SO_Y)=0$, we obtain the sequence
\begin{equation*}
 0 \to H^0(\SO_B(B)) \to \mathrm{Def}_{(Y,B)}(\mathbf k[\epsilon]) \xrightarrow{\mathrm{d}\nu} H^0(\SN_{Y/\PP^n}) \to 0.
\end{equation*}
Recall that $\mathrm{Def}_{\varphi}(\mathbf k[\epsilon])=H^0(\SN_{\varphi})$ and, from \cite[Lemma 3.3]{Gon}, there is a sequence
\begin{equation*}
 0 \to H^0(\SN_{\pi}) \to H^0(\SN_{\varphi}) \xrightarrow{\Psi_1 \oplus \Psi_2} H^0(\SN_{Y/\PP^N})\oplus H^0(\SN_{Y/\PP^N} \otimes \omega_Y (-1)).
\end{equation*}
Moreover, since we are assuming $n=2$, the restriction of $\pi$ to the ramification divisor $R$ is an isomorphism between $R$ and the branch divisor $B$. Then~\eqref{normalpi} implies that $\SN_{\pi}$ becomes $\SO_B(B)$ under this isomorphism and, as a consequence, there exist isomorphisms
\begin{equation}\label{cohom.normalpi.formula}
H^i(\SN_{\pi}) \simeq H^i(\SO_B(B))
\end{equation}
for all $i \geq 0$.

\smallskip
\noindent
Let $(\bar{Y}, \bar{B}) \in \mathrm{Def}_{(Y,B)}(\mathbf k[\epsilon])$ be a first order deformation of $(Y \subset \PP^N, B \in | \omega_Y^{-2}(2) |)$ and let $(\widetilde{X}, \widetilde{\varphi}) \in H^0(\SN_{\varphi})$ be the first order deformation of $(X, \varphi)$ associated to $(\bar{Y}, \bar{B})$ by $\mathrm{d}F$. From the construction we made for $F$, we see that $\widetilde{X}\xrightarrow{\widetilde{\varphi}} \PP^N \times \Delta$ factors through $\bar{Y} \hookrightarrow \PP^N \times \Delta$. Therefore $\mathrm{im}\,\widetilde{\varphi} =\bar{Y}$. Then, using \cite[Theorem 3.8 (2) and Propositions 3.11 and 3.12]{Gon}, we see that there is a commutative diagram
\begin{equation}\label{smart.diagram}
\xymatrix@C-10pt{
& 0 \ar[d] & 0 \ar[d] \\
0 \ar[r] & H^0(\SO_{B}(B)) \ar[d] \ar^\simeq[r] & H^0(\SN_{\pi}) \ar[d] \\
& \mathrm{Def}_{(Y,B)}(\mathbf k[\epsilon]) \ar[d]_{\mathrm{d}\nu} \ar[r]^{\mathrm{d}F} & H^0(\SN_{\varphi}) \ar[d]^{\Psi_1 \oplus \Psi_2} \\
0 \ar[r] & H^0(\SN_{Y/\PP^N}) \ar[d] \ar[r] & H^0(\SN_{Y/\PP^N})\oplus H^0(\SN_{Y/\PP^N} \otimes \omega_Y(-1)) \\
& 0&  .
}
\end{equation}
We see, from diagram \eqref{smart.diagram}, that if $\Psi_2=0$, then $\mathrm{d}F$ is an isomorphism.
\end{proof}
\begin{remark}\label{kernel.psi2}%**
%{\rm
Assuming the same hypotheses of Theorem \ref{Psi2=0} except hypothesis 
(6) (i.e., $\Psi_2$ need not be necessarily $0$), if  $H^1(\SO_B(B))=0$, then there exists an exact sequence
\begin{equation*}
 0 \to \mathrm{Def}_{(Y,B)}(\mathbf k[\epsilon]) \xrightarrow{\mathrm{d}F} H^0(\SN_{\varphi}) \xrightarrow{\Psi_2} H^0(\SN_{Y, \PP^N} \otimes \omega_Y(-1)) \to 0.
\end{equation*}
\end{remark}

\begin{proof}
This follows from the arguments in the proof of Theorem \ref{Psi2=0}.%}
\end{proof}

\begin{theorem}\label{Psi2<>0}
With the  notation of~\ref{setup} and of Proposition~\ref{morphism.miguel}, let $X$ be a smooth variety of general type of dimension $m \geq 2$ with ample and base--point--free canonical bundle, let $\varphi$ be its canonical morphism and assume that the degree of $\pi$ is $n=2$.
Assume furthermore that
\begin{enumerate}
\item $h^1(\mathcal O_Y)=h^{m-1}(\mathcal O_Y)=0$ (in particular, $Y$ is regular);
\item $h^1(\mathcal O_Y(1))=h^{m-1}(\mathcal O_Y(1))=0$;
\item the variety $X$ is unobstructed; and
\item $\Psi_2 \neq 0$.
\end{enumerate}
Then $\varphi$ is unobstructed and there exists a flat family of morphisms, $\Phi: \mathcal X \longrightarrow \mathbf P^N_T$ over $T$, where $T$
is a smooth algebraic curve with a distinguished point $0$,
such that
\begin{enumerate}
\item[(a)] $\Phi_0=\varphi$, and
\item[(b)] for any $t \in T, t \neq 0$, the morphism $\Phi_t: \mathcal X_t \longrightarrow \mathbf P^N$ is the canonical map of $\mathcal X_t$ and is finite and of degree $1$.
\end{enumerate}
Thus the canonical map of a variety corresponding to a general point of the component of $X$ in its moduli space is a finite morphism of degree $1$.
\end{theorem}

\begin{proof}
By hypothesis, $\omega_X$ is
base--point--free and by conditions (1) and (2), $H^1(\mathcal O_X)=H^{m-1}(\mathcal O_X)=0$. Then, since
$X$ is unobstructed, Corollary~\ref{canonical.remark}, (1) implies
that $\varphi$ is also unobstructed and Lemma~\ref{effective-1} implies that $\varphi$ has an algebraic formally universal deformation.  On the other hand, since $\Psi_2 \neq 0$, there exists $\nu \in H^0(\mathcal N_\varphi)$ such that $\Psi_2(\nu) \neq 0$. Since $n = 2$, this means that $\Psi_2(\nu)$ has rank $k > n/2-1$. Thus condition (1) and (2) of Theorem~\ref{onetoone2} are satisfied and the existence of $\Phi$ follows from Theorem~\ref{onetoone2} and Lemma~\ref{canonical.deforms.canonical}. Finally, since $X$ is unobstructed, the point of $X$ in its moduli belongs to a unique component of the moduli space.
The general point of this component has a canonical map which is a finite morphism of degree $1$, because such condition is open.
\end{proof}

\section{Canonical double covers of $\mathbf P^2$ blown--up at points in general position}\label{construct.from.blownup.section}

In this section we will apply the results of Sections~\ref{intrinsic.section} and~\ref{deformations.canonical.general.section} to study the deformations of canonical double covers of non--minimal rational surfaces obtained by blowing--up $\mathbf P^2$ at points in general position. We will need to introduce some further notation:

\begin{noname}\label{notation.blowups}
{\bf Notation.} {\rm In addition to the notation introduced in~\ref{setup}, throughout the remaining of the article
\begin{enumerate}
 \item we will denote by $Y$ a smooth rational surface obtained by blowing up $\mathbf P^2$ at a set of points $S=\{x_1, \dots, x_s\}$ ($s \geq 1$), in sufficiently general position and
$p:  Y \longrightarrow \mathbf P^2$ will be this blowing--up;
\item we will denote by $\mathfrak m_{x_i}$ the ideal sheaf of $x_i$ and $\mathfrak m=\mathfrak m_{x_1} \otimes \cdots \otimes \mathfrak m_{x_s}$;
\item we will call $L$ the pullback by $p$ of a line in $\mathbf P^2$, $K_Y$ the canonical divisor of $Y$, $E_1, \dots, E_s$ the exceptional divisors of $p$ and $E=E_1+ \cdots + E_s$.
\end{enumerate}}
\end{noname}

We will focus on the canonical covers of embeddings of $Y$ by complete linear series of the form
$|dL -E_1 - \cdots -E_s|=|dL-E|$, i.e., the linear series corresponding to the system of curves in $\mathbf P^2$ of degree $d$, passing through the points of $S$. Thus we need to know for which values of $d$ and $s$ we can assure $|dL -E|$ to be very ample (for this to happen necessarily $d \geq 2$). Next, among those values of $d$ and $s$, we want to know for which values there exists a smooth canonical double cover of the image of $Y$ by the embedding induced by $|dL -E|$. This we do in the next Propositions~\ref{blowup.embeddings} and~\ref{exist.emb.cover} and Corollary~\ref{exist.emb.cover.cor}:

\begin{proposition}\label{blowup.embeddings}
Let $L$, $s$, $d$ and $E_i$ be as in~\ref{notation.blowups}. Then $|dL -E|$ is very ample if and only if one of the following occurs:
\begin{enumerate}
\item[(a)] $d=2$ and $s =1$, or
\item[(b)] $d=3$ and $s \leq 6$, or
\item[(c)] $d=4$ and $s \leq 10$, or
\item[(d)] $d \geq 5$  and $s \leq \frac{1}{2}d^2+\frac{3}{2}d-5$.
\end{enumerate}
\end{proposition}

\begin{proof}
First let us assume $d=2$. In this case $|dL -E|$ is very ample if $s=1$, since in this case $Y$ is the Hirzebruch surface $\mathbf F_1$ and $|dL -E|$ embeds $\mathbf F_1$ as a cubic scroll in $\mathbf P^4$. On the other hand, if $d=2$ and $s=2$, $|dL -E|$ cannot be very ample, for otherwise it would embed $Y$ as a quadric surface in $\mathbf P^3$. Likewise, if $d=2$, $|dL -E|$ cannot be very ample if $s \geq 3$, for in this case the dimension of the linear system is $2$ or less. This completes the proof of (a).

\smallskip
\noindent
Now let us assume $d=3$. Then $|dL -E|$ is very ample if and only if $s \leq 6$ (see \cite[Theorem V.4.6]{Hart}).
This proves (b).

\smallskip
\noindent
Now let us assume $d=4$. In this case, if $s=10$, then $|dL -E|$ is also very ample (see~\cite{Ionescu}). This implies that $|dL -E|$ is very ample for $d=4$ and $s \leq 9$. On the other hand, if $s \geq 11$, $|dL -E|$ cannot be very ample. Indeed, since $\{x_1,\dots,x_s\}$ are in general position, they impose independent conditions on quartics, so $|dL -E|$ could be very ample only if $s=11$. However, if that happened, $Y$ will be isomorphic to a smooth quintic in $\mathbf P^3$, which is a surface of general type. This completes the proof of (c).

\smallskip
\noindent
 Finally let us assume $d \geq 5$. If $\frac{(d+2)(d+1)}{2}-s \geq 6$, then $|dL -E|$ is very ample by~\cite[Theorem 2.3]{dAH}. This inequality is equivalent to
$s \leq \frac{1}{2}d^2+\frac{3}{2}d-5$.
On the other hand, if $d \geq 5$ and $s \geq \frac{1}{2}d^2+\frac{3}{2}d-4$, we argue as for $d=4$. Since $\{x_1,\dots,x_s\}$ impose independent conditions on curves of degree $d$, $|dL -E|$ could be very ample only if $s = \frac{1}{2}d^2+\frac{3}{2}d-4$ or $s = \frac{1}{2}d^2+\frac{3}{2}d-3$. In the case $s = \frac{1}{2}d^2+\frac{3}{2}d-3$, $|dL -E|$ would embed $Y$ in $\mathbf P^3$ with degree $d^2-s= \frac{1}{2}d^2-\frac{3}{2}d+3 \geq 8$, so again $Y$ should be a surface of general type. In the case $s = \frac{1}{2}d^2+\frac{3}{2}d-4$, $|dL -E|$ would embed $Y$ in $\mathbf P^4$ with degree $d^2-s= \frac{1}{2}d^2-\frac{3}{2}d+4 \geq 9$. By \cite[Proposition 4.3]{Alex}, this might only be possible if $d^2-s=\frac{1}{2}d^2-\frac{3}{2}d+4=9$ and $K_Y^2=-1$. But if $d^2-s=\frac{1}{2}d^2-\frac{3}{2}d+4=9$, then $d=5$ and $s=16$, so $K_Y^2=-7$ and we get a contradiction. This completes the proof of (d). \end{proof}

\noindent
Since a canonical cover of $Y$ is branched along $|\omega_Y^{-2}(2)|$ (see~\eqref{branch.of.canonical.double.cover}) and we want to see for what values of $d$ and $s$, the surface $Y$ admits a smooth canonical double cover, we will study when $|\omega_Y^{-2}(2)|$ has smooth divisors. This we do in Proposition~\ref{exist.emb.cover}, in which we will use the following lemma, whose proof was communicated to us by Brian Harbourne.

\begin{lemma}\label{Harbourne}
If $d=5$ and $s=14$,
then
$|-2K_Y+2dL-2E|$ is base--point--free.
\end{lemma}

\begin{proof}
Observe that
\begin{equation*}\label{canon.hyper2.1}
-2K_Y+2dL-2E=(2d+6)L-4E.
\end{equation*}
So for $(d,s)=(5,14)$ we have to prove that the linear system $|16 L -4E_1-\cdots-4E_{14}|$ is base--point--free.
Through the $14$ general points in $\mathbf P^2$ there is a unique smooth quartic $C'$. Let $C \in |4L -E_1 -\cdots -E_{14}|$ be its proper transform on the blow up; then $|4C|=|16 L -4E_1-\cdots-4E_{14}|$.

\noindent
Consider the sequence
\begin{equation*}
\xymatrix@C-5pt{0 \ar[r]& \SO_Y(3C)\ar[r] & \SO_Y(4C) \ar[r] & \SO_C(4C) \ar[r] & 0.}
\end{equation*}
It is known (\cite[Lemma 7.1 and Theorem 8.1]{CilibertoMiranda2}) that the linear system $|12L-3E_1-\cdots -3E_{14}|$ is non--special, so $h^1(\SO_Y(3C))=0$ and $h^0(\SO_Y(3C))=7$.
Note that $\mathrm{deg}(\SO_C(4C))=8=2g(C)+2$, hence $|H^0(\SO_C(4C))|$ is base--point--free, $h^1(\SO_C(4C))=0$ and $h^0(\SO_C(4C))=6$. Therefore $H^0(\SO_Y(4C)) \to H^0(\SO_C(4C))$ is surjective, $h^1(\SO_Y(4C))=0$ (i.e. $|16L-4E_1-\cdots -4E_{14}|$ is non--special) and $h^0(\SO_Y(4C))=13$. Since $H^0(\SO_Y(4C)) \to H^0(\SO_C(4C))$ is surjective and $|H^0(\SO_C(4C))|$ is base--point--free, $|4C|$ is also base--point--free.
\end{proof}

\begin{proposition}\label{exist.emb.cover}
Let $Y$, $L$, $s$, $d$ and $E_i$ be as in~\ref{notation.blowups}.
\begin{enumerate}
\item When $2 \leq d \leq 4$, the linear system $|dL -E|$ is very ample and
the linear system $|-2K_Y+2dL-2E|$ has a smooth member if and only if
\begin{enumerate}
\item[(a)] $d=2$ and $s =1$, or
\item[(b)] $d=3$ and $s \leq 6$, or
\item[(c)] $d=4$ and $s \leq 10$.
\end{enumerate}
\item When $d \geq 5$,
\begin{enumerate}
\item[(d)] $d=5$ and $s \leq 14$, and
\item[(e)] $d \geq 6$ and $s \leq \frac{1}{5}d^2+\frac{13}{10}d+\frac{21}{10}$
\end{enumerate}
\item[] are sufficient conditions and,
\begin{enumerate}
\item[(f)] $s < \frac{1}{5}d^2+\frac{3}{2}d+\frac{14}{5}$
\end{enumerate}
\item[] is a necessary condition, for $|dL -E|$ to be very ample and
$|-2K_Y+2dL-2E|$ to have a smooth member.
\end{enumerate}
\end{proposition}

\begin{proof}
We prove all sufficient conditions first. They would follow from 
Proposition~\ref{blowup.embeddings} if we also prove that $|-2K_Y+2dL-2E|$ has a smooth member when (a), (b), (c), (d) or (e) are satisfied.  To see this, recall that
\begin{equation}\label{canon.hyper2}
-2K_Y+2dL-2E=(2d+6)L-4E.
\end{equation}
By Bertini's Theorem, there exists a smooth member in $|(2d+6)L-4E|$ if this linear system is base--point--free.
We argue first when $(d,s) \neq (5,14)$.
The fact that $|(2d+6)L-4E|$ is base--point--free will follow from Castelnuovo--Mumford regularity (see \cite{Mu}), if we see that the line bundle $\SO_Y(-4E)$ is $(2d+6)$--regular with respect to $L$. So, by definition, we have to prove that $H^1(\SO_Y((2d+5)L-4E))=0$ and $H^2(\SO_Y((2d+4)L-4E))=0$. The latter vanishing is the equivalent to $H^2(\SO_{\PP^2}(2d+4) \otimes \mathfrak m^4)=0$. The sequence
\begin{equation*}
H^1(\SO_{\PP^2}(2d+4) \mid_{\Sigma}) \to H^2(\SO_{\PP^2}(2d+4) \otimes \mathfrak m^4)) \to H^2(\SO_{\PP^2}(2d+4)),
\end{equation*}
where $\Sigma$ is the third infinitesimal neighborhood of $S$, implies $H^2(\SO_{\PP^2}(2d+4) \otimes \mathfrak m^4)=0$. For $H^1(\SO_Y((2d+5)L-4E))=0$, we use~\cite[Theorems 2.4 and 5.2]{CilibertoMiranda}. Since $2d+5 >9$, and, by \cite[Theorem 2.4]{CilibertoMiranda}, there do not exist homogeneous, $(-1)$--special systems of multiplicity $4$ and degree bigger than $9$, it follows, from \cite[Theorem 5.2]{CilibertoMiranda}, that $H^1(\SO_Y((2d+5)L-4E))=0$ is equivalent to $\frac{(2d+7)(2d+6)}{2}-10s \geq 0$. The latter inequality is equivalent to $s \leq \frac{1}{5}d^2+\frac{13}{10}d+\frac{21}{10}$.
Note also that for $d=3$, the inequality $s \leq \frac{1}{5}d^2+\frac{13}{10}d+\frac{21}{10}$ becomes $s \leq 7$, for $d=4$, it becomes $s \leq 10$ and for $d=5$, it becomes $s \leq 13$. This proves the existence of a smooth member in $|-2K_Y+2dL-2E|$ if $(d,s)$ satisfies (a), (b), (c), (d)  or (e), with the exception of $(d,s)=(5,14)$.

\smallskip
\noindent
The existence of a smooth member in $|-2K_Y+2dL-2E|$ when $(d,s)=(5,14)$
follows from Lemma~\ref{Harbourne}.
This completes the proof of the sufficient conditions stated in the proposition.

\smallskip
\noindent
Now we prove the necessary conditions stated in the proposition. If $2 \leq d \leq 4$, they follow from Proposition~\ref{blowup.embeddings}.  Now let $d \geq 5$. The inequality $s \geq \frac{1}{5}d^2+\frac{3}{2}d+\frac{14}{5}$ is equivalent to $\frac{(2d+8)(2d+7)}{2}-10s \leq 0$. In particular, $s \geq 16$. Then~\cite[Theorems 2.4 and 5.2]{CilibertoMiranda} imply $|(2d+6)L-4E|$ is empty.
\end{proof}

\begin{corollary}\label{exist.emb.cover.cor}
\begin{enumerate}
 \item When $2 \leq d \leq 4$, there exist smooth surfaces of general type $X$ with ample and base--point--free canonical bundle and whose canonical morphism $\varphi$ maps $2:1$ onto a surface $Y$, embedded in projective space by the morphism $i$ induced by the linear system $|dL -E|$ if and only if $s$ and $d$ are as in (a), (b) or (c) of Proposition~\ref{exist.emb.cover}.
\item When $d \geq 5$, there exist surfaces $X$ with the same properties described in (1) above if $s$ and $d$ are as in (d) or (e) of Proposition~\ref{exist.emb.cover} and there do not exist if $s$ and $d$  are as in (f) of Proposition~\ref{exist.emb.cover}.
\end{enumerate}
\end{corollary}

\begin{proof}
First we study when the surfaces $X$ described in the statement do exist. Assume that $s$ and $d$ are as in (a), (b), (c), (d)  or (e) of Proposition~\ref{exist.emb.cover}. Let $B$ be a smooth curve in $|-2K_Y+2dL-2E|$. Since on $Y$ numerical and linear equivalence are the same, there is only one line bundle, which is $\omega_Y^{-1}(1)$, whose square is $\mathcal O_Y(B)$ (recall that $\mathcal O_Y(1)=\mathcal O_Y(dL -E)$). Thus the double cover $\pi: X \longrightarrow Y$, of $Y$ branched along $B$, satisfies $\omega_X=\pi^*\mathcal O_Y(1)$, so $\omega_X$ is therefore ample and base--point--free. Moreover, since $p_g(Y)=0$, the canonical morphism $\varphi$ of $X$ factors as $\varphi= i \circ \pi$. Finally, since $B$ is smooth, so is $X$.

\smallskip
\noindent Now we find the conditions under which the surfaces $X$ in the statement do not exist. If  $2 \leq d \leq 4$, assume that $s$ and $d$ are not as in (a), (b) or (c) of  Proposition~\ref{exist.emb.cover}; if $d \geq 5$, assume that $s$ and $d$ satisfy (f) of  Proposition~\ref{exist.emb.cover}. Suppose furthermore that there exists a smooth surface $X$ as described in the statement of the corollary. Then we would arrive to a contradiction, for either the image $Y$ of $\varphi$ cannot be embedded by $|dL-E|$, or the fact that $K_X$ is the pullback of $dL-E$ by $\pi$ would imply that the branch locus of $\pi$ is a smooth curve in $|(2d+6)L-4E|$, which in particular would not be empty.
\end{proof}

\begin{remark}\label{branch.of.canonical.double.cover.remark}
{\rm By~\eqref{branch.of.canonical.double.cover}, a smooth surface of general type $X$ with ample and base--point--free canonical bundle and whose canonical morphism $\varphi$ maps $2:1$ onto a surface $Y$ embedded in projective space has necessarily a smooth member of $|\omega_Y^{-2}(2)|$ as its branch locus; therefore any smooth surface of general type $X$ with ample and base--point--free canonical bundle and whose canonical morphism $\varphi$ maps $2:1$ onto a surface $Y$, embedded in projective space by the linear system $|dL -E|$ with $s$ and $d$ as in (a), (b), (c), (d) or (e) of Proposition~\ref{exist.emb.cover}, is constructed as in the proof of Corollary~\ref{exist.emb.cover.cor}.}
\end{remark}

\noindent
Now that
we have constructed canonical double covers $\varphi$ of $Y$ our next task will be to decide when the homomorphism $\Psi_2$ associated to $\varphi$ (recall Proposition~\ref{morphism.miguel}) is different from $0$ and when is $0$. This way, using Theorems~\ref{Psi2=0} and~\ref{Psi2<>0}, we may draw conclusions about the deformations of the covers exhibited in Corollary~\ref{exist.emb.cover.cor}. The homomorphism $\Psi_2$ goes to Hom$(\mathcal I/\mathcal I^2,\omega_Y(-1))$. We will see, as a consequence of Lemma~\ref{isom.Hom.Ext}, that, in our setting, Ext$^1(\Omega_Y,\omega_Y(-1))$ is isomorphic to Hom$(\mathcal I/\mathcal I^2,\omega_Y(-1))$, so we will study the former in order to handle the latter. We will accomplish this in Corollary~\ref{exist.abstr.carpets} and Theorem~\ref{exist.embedd.carpets}, where we will also obtain a result on the existence of ``general type'' double structures of dimension $2$ and their ``canonical'' morphisms.
\medskip

\noindent As we will see in Corollary~\ref{exist.abstr.carpets},
Ext$^1(\Omega_Y,\omega_Y(-1))$ is related to a multiplication map of global sections on $\mathbf P^2$. Therefore we continue with this lemma:

\begin{lemma}\label{surjectivity} Let $d \geq 2$. Let $S=\{x_1,\dots,x_s\}$ be as in~\ref{notation.blowups}. Let $\alpha$ be the multiplication map of global sections on $\mathbf P^2$
\begin{equation*}
 H^0(\mathcal O_{\mathbf P^2}(d-1) \otimes \mathfrak m) \otimes H^0(\mathcal O_{\mathbf P^2}(1)) \overset{\alpha} \longrightarrow H^0(\mathcal O_{\mathbf P^2}(d) \otimes \mathfrak m).
\end{equation*}
\begin{enumerate}
\item When $2 \leq d \leq 4$, $\alpha$ surjects if and only if
 \begin{enumerate}
  \item[(a)] $d=2$ and $s=1$ or $s \geq 6$; or
  \item[(b)] $d=3$ and $s \leq 4$ or $s \geq 10$; or
  \item[(c)] $d=4$ and $s \leq 7$ or $s \geq 15$.
 \end{enumerate}
\item When $d \geq 5$, 
  \begin{enumerate}
   \item[(d)] if $s \leq \frac{1}{2}d^2-\frac{1}{2}d+1$ or $s \geq \frac{1}{2}d^2+\frac{3}{2}d+1$, then $\alpha$ surjects; and
   \item[(e)] if $\frac{1}{2}d^2-\frac{1}{2} < s <  \frac{1}{2}d^2+\frac{3}{2}d+1$, then $\alpha$ does not surject.
  \end{enumerate}
\end{enumerate}
\end{lemma}

\begin{proof}
First we prove the result when $d=2$. If $d=2$ and $s=1$, $H^1(\mathcal O_{\mathbf P^2}(d-2) \otimes \mathfrak m)=H^1(\mathcal O_{\mathbf P^2} \otimes \mathfrak m)=0$. Indeed, $S$ consists of just one point $x_1$, and
the evaluation map
\begin{equation*}
 H^0(\mathcal O_{\mathbf P^2}) \longrightarrow H^0(\mathcal O_{\mathbf P^2}|_{x_1})
\end{equation*}
is surjective, so the claim follows.
Moreover, the exact sequence
\begin{equation*}
 H^1(\mathcal O_{\mathbf P^2}(-1)|_{x_1}) \longrightarrow H^2(\mathcal O_{\mathbf P^2}(-1) \otimes \mathfrak m) \longrightarrow H^2(\mathcal O_{\mathbf P^2}(-1))
\end{equation*}
implies that $H^2(\mathcal O_{\mathbf P^2}(-1) \otimes \mathfrak m)=0$. Then~\cite[p.41, Theorem 2]{Mu} implies the surjectivity of $\alpha$.
On the other hand, if $d=2$ and $s \geq 6$, $H^0(\mathcal O_{\mathbf P^2}(d) \otimes \mathfrak m)=0$ since the points of $S$ are in general position, so $\alpha$ surjects trivially.

\smallskip

\noindent Now let $d=2$ and $2 \leq s \leq 5$. If we compare the dimensions of $H^0(\mathcal O_{\mathbf P^2}(d-1) \otimes \mathfrak m) \otimes H^0(\mathcal O_{\mathbf P^2}(1))$ and $H^0(\mathcal O_{\mathbf P^2}(d) \otimes \mathfrak m)$ we easily see that the dimension of the former is strictly smaller than the dimension of the latter so obviously $\alpha$ cannot be surjective.

\smallskip
\noindent 
Now we deal with the remaining of (1) and with (2). What we have to prove is equivalent to showing that 
\begin{enumerate}
\item[] if $d \geq 3$ and $s \leq \frac{1}{2}d^2-\frac{1}{2}d+1$ or $s \geq \frac{1}{2}d^2+\frac{3}{2}d+1$, then $\alpha$ surjects; and
\item[] if $d \geq 3$ and $\frac{1}{2}d^2-\frac{1}{2} < s <  \frac{1}{2}d^2+\frac{3}{2}d+1$, then $\alpha$ does not surject. 
\end{enumerate}

\smallskip
\noindent
First we see what happens if $d \geq 3$ and $s \leq \frac{1}{2}d^2-\frac{1}{2}d$. In this case $H^1(\mathcal O_{\mathbf P^2}(d-2) \otimes \mathfrak m)=0$. Indeed, the points of $S$ are in general position, so they impose independent conditions on curves of degree $d-2$. Since $s \leq \frac{d(d-1)}{2}$ this means that the evaluation map
\begin{equation*}
 H^0(\mathcal O_{\mathbf P^2}(d-2)) \longrightarrow H^0(\mathcal O_{\mathbf P^2}(d-2)|_S)
\end{equation*}
is surjective, so the claim follows.
On the other hand the exact sequence
\begin{equation*}
 H^1(\mathcal O_{\mathbf P^2}(d-3)|_S) \longrightarrow H^2(\mathcal O_{\mathbf P^2}(d-3) \otimes \mathfrak m) \longrightarrow H^2(\mathcal O_{\mathbf P^2}(d-3))
\end{equation*}
implies that $H^2(\mathcal O_{\mathbf P^2}(d-3) \otimes \mathfrak m)=0$. Then~\cite[p.41, Theorem 2]{Mu} implies the surjectivity of $\alpha$.

\smallskip
\noindent
Now we study the case in which $d \geq 3$ and $s = \frac{1}{2}d^2-\frac{1}{2}d+1$. We claim that the linear system $|(d-1)L-E|$
on $Y$ is base--point--free. Indeed, if $d=3$, then $s=4$, so $|(d-1)L-E|$ corresponds to a linear system of conics passing through $4$ general points of $\mathbf P^2$. Such a system has no unassigned base points, so $|(d-1)L-E|$ is base--point--free. Now, if $d \geq 4$, since $\frac{(d+1)d}{2}-s \geq d-1\geq 3$ and the points are taken in general position, ~\cite[3.3]{Cop} implies the claim.
Now let $l$ be a straight line in $\mathbf P^2$ not passing through any point of $S$. Then the linear system of curves of $\mathbf P^2$ of degree $d-1$ passing through $\{x_1,\dots,x_s\}$ restricts to a base--point--free linear system $|W_{d-1}|$ on $l$ of divisors of degree $d-1$. Let $V_{d-1}= H^0(\mathcal O_{\mathbf P^2}(d-1) \otimes \mathfrak m)$ and $V_d=H^0(\mathcal O_{\mathbf P^2}(d) \otimes \mathfrak m)$.
Consider the following commutative diagram:
\begin{equation}\label{restrict.tocurves}
\xymatrix{0 \ar[r] &
V_{d-1} \otimes H^0(\mathcal O_{\mathbf P^2}) \ar[r] \ar[d] &  V_{d-1} \otimes H^0(\mathcal O_{\mathbf P^2}(1)) \ar[r] \ar[d]^\alpha & V_{d-1} \otimes H^0(\mathcal O_{l}(1)) \ar[r]
\ar[d]^{\alpha'} & 0
\\
0 \ar[r] &
V_{d-1} \ar[r] &  V_d \ar[r] & H^0(\mathcal O_{l}(d))
,}
\end{equation}
where the vertical arrows are the obvious multiplication maps of global sections.
The left--hand--side vertical arrow is evidently an isomorphism. Thus, if we prove that $\alpha'$ surjects, we will prove that $\alpha$ also surjects. Now $\alpha'$ is the composition of these two maps:
\begin{eqnarray*}
 V_{d-1} \otimes H^0(\mathcal O_{l}(1)) & \overset{\alpha_1'}\longrightarrow &  W_{d-1} \otimes H^0(\mathcal O_{l}(1)) \cr
W_{d-1} \otimes H^0(\mathcal O_{l}(1)) & \overset{\alpha_2'} \longrightarrow & H^0(\mathcal O_{l}(d)),
\end{eqnarray*}
where $\alpha'_1$ surjects. Thus we just need to show that $\alpha_2'$ also surjects. For that we define the vector bundles $M_1$ and $M_2$ as the kernels of the following evaluation maps of global sections (note that the evaluation maps are both surjective because $\mathcal O_{l}(d-1)$ is globally generated and $|W_{d-1}|$ is base--point--free):

\begin{equation}\label{M1}
0 \longrightarrow
M_1
\longrightarrow
 H^0(\mathcal O_{l}(d-1)) \otimes \mathcal O_{l}  \longrightarrow  \mathcal O_{l}(d-1)  \longrightarrow  0
\end{equation}
\begin{equation} \label{M2}
0 \longrightarrow
M_2 \longrightarrow
W_{d-1} \otimes \mathcal O_{l}  \longrightarrow  \mathcal O_{l}(d-1)  \longrightarrow  0.
\end{equation}

If we twist~\eqref{M2} with $\mathcal O_{l}(1)$ and take cohomology, we see that the surjectivity of $\alpha'_2$ is equivalent to the vanishing of $H^1(M_2 \otimes \mathcal O_{l}(1))$.
To prove this vanishing we observe that the exact sequences~\eqref{M1} and~\eqref{M2} fit in the following commutative diagram:
\begin{equation*}
\xymatrix{
          & 0 \ar[d]          & 0 \ar[d]          &        0 \ar[d]                                                           &   \\
0 \ar[r]  & M_2 \ar[r] \ar[d] & M_1 \ar[r] \ar[d] & H^0(\mathcal O_l(d-1))/W_{d-1} \otimes \mathcal O_l \ar[r] \ar[d] & 0 \\
0 \ar[r]  & W_{d-1} \otimes \mathcal O_l \ar[r] \ar[d] & H^0(\mathcal O_l(d-1)) \otimes \mathcal O_l \ar[r] \ar[d] & H^0(\mathcal O_l(d-1))/W_{d-1} \otimes \mathcal O_l \ar[r] \ar[d] & 0 \\
0 \ar[r]  & \mathcal O_l(d-1) \ar[r] \ar[d] & \mathcal O_l(d-1) \ar[r] \ar[d] & 0\\
          &  0                              &  0                             &.}
\end{equation*}
Since $h^0(\mathcal O_{\mathbf P^2}(d-2)) < s=\frac{1}{2}d^2-\frac{1}{2}d+1$
and $\{x_1,\dots,x_s\}$ are in general position, $H^0(\mathcal O_{\mathbf P^2}(d-2) \otimes \mathfrak m)=0$, and $V_{d-1}$ is isomorphic to $W_{d-1}$ which has therefore dimension $d-1$. On the other hand, taking cohomology at~\eqref{M1} we see that $H^0(M_1)=H^1(M_1)=0$. Since a vector bundle on $l \simeq \mathbf P^1$ splits as a direct sum of line bundles, then $M_1=\mathcal O_{l}^{\oplus d-1}(-1)$. Thus,
from the top horizontal exact sequence we get that $M_2$ fits in
\begin{equation*}\label{M1M2}
 0 \longrightarrow M_2 \longrightarrow \mathcal O_{l}^{\oplus d-1}(-1) \longrightarrow \mathcal O_{l} \longrightarrow 0.
\end{equation*}
Then $M_2 \simeq \mathcal O_{l}^{\oplus d-3}(-1) \oplus \mathcal O_{l}(-2)$. Then $H^1(M_2 \otimes \mathcal O_{l}(1))=0$ as wished.

\smallskip
\noindent
Now we deal with the range $d \geq 3$ and $\frac{1}{2}d^2-\frac{1}{2} < s <  \frac{1}{2}d^2+\frac{3}{2}d+1$.
In this case the dimension of $H^0(\mathcal O_{\mathbf P^2}(d-1) \otimes \mathfrak m) \otimes H^0(\mathcal O_{\mathbf P^2}(1))$ is smaller than the dimension of $H^0(\mathcal O_{\mathbf P^2}(d) \otimes \mathfrak m)$, so obviously $\alpha$ cannot surject.

\smallskip
\noindent
Finally we study the case in which $d \geq 3$ and $s \geq \frac{1}{2}d^2+\frac{3}{2}d+1$. In this situation $H^0(\mathcal O_{\mathbf P^2}(d) \otimes \mathfrak m) =0$ because the points of $S$ are in general position, so $\alpha$ obviously surjects.
\end{proof}

\begin{corollary}\label{exist.abstr.carpets}
Let $Y$, $d$, $s$, $E_i$ and $L$ be as in~\ref{notation.blowups}, let $d \geq 2$ and let $ M=\mathcal O_Y(dL-E_1-\cdots-E_s)$.

\begin{enumerate}
\item When $2 \leq d \leq 4$, $\mathrm{Ext}^1(\Omega_Y,\omega_Y \otimes M^\vee) \neq 0$ (thus, there exist double structures on $Y$ with conormal bundle $\omega_Y \otimes M^\vee$) if and only if
 \begin{enumerate}
  \item[(a)] $d=2$ and $s \geq 2$; or
  \item[(b)] $d=3$ and $s \geq 5$; or
  \item[(c)] $d=4$ and $s \geq 8$.
 \end{enumerate}
\item When $d \geq 5$,
  \begin{enumerate}
   \item[(d)]  if $s \leq \frac{1}{2}d^2-\frac{1}{2}d+1$, then $\mathrm{Ext}^1(\Omega_Y,\omega_Y \otimes M^\vee) = 0$ (thus, there do not exist double structures on $Y$ with conormal bundle $\omega_Y \otimes M^\vee$); and
   \item[(e)] if $s > \frac{1}{2}d^2-\frac{1}{2}$, then $\mathrm{Ext}^1(\Omega_Y,\omega_Y \otimes M^\vee) \neq 0$ (thus,
there exist double structures on $Y$ with conormal bundle $\omega_Y \otimes M^\vee$).
  \end{enumerate}
\end{enumerate}
\end{corollary}

\begin{proof}
The result follows from Lemma~\ref{surjectivity} once we see that the cokernel of the multiplication map $\alpha$ injects in Ext$^1(\Omega_Y,\omega_Y  \otimes M^\vee)$. By duality, Ext$^1(\Omega_Y,\omega_Y  \otimes M^\vee)$ is isomorphic to $H^1(\Omega_Y \otimes M)^{\vee}$, so we study the latter. For this we use the sequence
\begin{equation}\label{rel.cotangent.seq}
0 \longrightarrow p^*\Omega_{\mathbf P^2} \longrightarrow \Omega_Y \longrightarrow \Omega_{Y/\mathbf P^2} \longrightarrow 0,
\end{equation}
tensored by $M$.
Recall that
\begin{equation}\label{rel.cotangent}
\Omega_{Y/\mathbf P^2}=\mathcal O_{E_1}(2E_1) \oplus \cdots \oplus \mathcal O_{E_s}(2E_s)
\end{equation}
Then $\Omega_{Y/\mathbf P^2} \otimes M= \mathcal O_{\mathbf P^1}^{\oplus s}(-1)$, so $H^0(\Omega_{Y/\mathbf P^2}(1))=H^1(\Omega_{Y/\mathbf P^2}(1))=0$, so
\begin{equation*}
H^1(\Omega_Y \otimes M) = H^1(p^*\Omega_{\mathbf P^2}\otimes M).
\end{equation*}
Now to compute $H^1(p^*\Omega_{\mathbf P^2}\otimes M)$ we
lift the Euler sequence from $\mathbf P^2$ to $Y$, twist with $M$, take cohomology and get
\begin{multline}\label{lifted.Euler}
0 \to H^0(p^*\Omega_{\mathbf P^2}\otimes M) \to H^0(p^*\mathcal O_{\mathbf P^2}(1)) \otimes H^0(p^*\mathcal O_{\mathbf P^2}(d-1) \otimes \mathcal O_Y(-E)) \overset{\tilde \alpha} \to H^0(M) \\
\to H^1(p^*\Omega_{\mathbf P^2}\otimes M)   \to H^0(p^*\mathcal O_{\mathbf P^2}(1)) \otimes H^1(p^*\mathcal O_{\mathbf P^2}(d-1) \otimes \mathcal O_Y(-E)).
\end{multline}
 Using projection formula we can compute the global sections in~\eqref{lifted.Euler} by pushing down the vector bundles back to $\mathbf P^2$, so we can identify $\tilde \alpha$ with the map $\alpha$.

\smallskip
\noindent Now we deal with 
the different cases in the statement. Note that all we need to prove is, on the one hand, the statement regarding $d=2$ and, on the other hand, the following:
\begin{enumerate}
\item[] if $d \geq 3$ and $s \leq \frac{1}{2}d^2-\frac{1}{2}d+1$, then $\mathrm{Ext}^1(\Omega_Y,\omega_Y \otimes M^\vee) = 0$; and
\item[] if $d \geq 3$ and $s > \frac{1}{2}d^2-\frac{1}{2}$, then $\mathrm{Ext}^1(\Omega_Y,\omega_Y \otimes M^\vee) \neq 0$.
\end{enumerate}

\smallskip
\noindent
If $d=2$ and $s=1$, Lemma~\ref{surjectivity} says that $\alpha$  surjects. On the other hand, in this case
\begin{equation*}
 H^1(p^*\mathcal O_{\mathbf P^2}(d-1) \otimes \mathcal O_Y(-E))= H^1(\mathcal O_{\mathbf P^2}(d-1) \otimes \mathfrak m)=H^1(\mathcal O_{\mathbf P^2}(1) \otimes \mathfrak m_{x_1})=0,
\end{equation*}
 because the evaluation map
\begin{equation*}
 H^0(\mathcal O_{\mathbf P^2}(1)) \longrightarrow H^0(\mathcal O_{\mathbf P^2}(1)|_{x_1})
\end{equation*}
obviously surjects.
Thus $H^1(p^*\Omega_{\mathbf P^2}\otimes M)$ and therefore Ext$^1(\Omega_Y,\omega_Y  \otimes M^\vee)$, are zero.

\smallskip
\noindent
If $d \geq 3$ and $s \leq \frac{1}{2}d^2-\frac{1}{2}d+1$ Lemma~\ref{surjectivity} says that $\alpha$  surjects. On the other hand, $H^1(p^*\mathcal O_{\mathbf P^2}(d-1) \otimes \mathcal O_Y(-E))= H^1(\mathcal O_{\mathbf P^2}(d-1) \otimes \mathfrak m)$ by the projection formula and the Leray spectral sequence.
The latter cohomology group vanishes because $x_1,\dots, x_s$ are in general position and $s \leq \frac{1}{2}d^2+\frac{1}{2}d$.  Thus $H^1(p^*\Omega_{\mathbf P^2}\otimes M)$ and therefore Ext$^1(\Omega_Y,\omega_Y  \otimes M^\vee)$, are zero.

\smallskip
\noindent If $d=2$ and $2 \leq s \leq 5$ or if $d\geq 3$ and $\frac{1}{2}d^2-\frac{1}{2} < s <  \frac{1}{2}d^2+\frac{3}{2}d+1$, Lemma~\ref{surjectivity} says that $\alpha$ does not surject, so $H^1(p^*\Omega_{\mathbf P^2}\otimes M)$, and therefore Ext$^1(\Omega_Y,\omega_Y  \otimes M^\vee)$, are not zero.

\smallskip
\noindent
Finally, if $d=2$ and $s \geq 6$ or if $d \geq 3$ and
 $s \geq \frac{1}{2}d^2+\frac{3}{2}d+1$, Lemma~\ref{surjectivity} says that $\alpha$ surjects. Then, using again the projection formula and the Leray spectral sequence we can identify $H^1(p^*\Omega_{\mathbf P^2}\otimes M)$ with the kernel of the map
\begin{equation*}
 H^1(\mathcal O_{\mathbf P^2}(d-1) \otimes \mathfrak m) \otimes H^0(\mathcal O_{\mathbf P^2}(1))  \longrightarrow H^1(\mathcal O_{\mathbf P^2}(d) \otimes \mathfrak m).
\end{equation*}
Now, since in this case both $H^0(\mathcal O_{\mathbf P^2}(d-1) \otimes \mathfrak m)$ and  $H^0(\mathcal O_{\mathbf P^2}(d) \otimes \mathfrak m)$ vanish, the dimension of $H^1(\mathcal O_{\mathbf P^2}(d-1) \otimes \mathfrak m) \otimes H^0(\mathcal O_{\mathbf P^2}(1))$ is $3(s-\frac{(d+1)d}{2})$, which is larger than $s-\frac{(d+2)(d+1)}{2}$, which is the dimension of $H^1(\mathcal O_{\mathbf P^2}(d) \otimes \mathfrak m)$. Thus $H^1(p^*\Omega_{\mathbf P^2}\otimes M)$, and therefore Ext$^1(\Omega_Y,\omega_Y  \otimes M^\vee)$, are not zero.

\smallskip

\noindent Finally, the remarks about double structures follow from~\cite[Theorem 1.2]{BE}.
\end{proof}

\begin{lemma}\label{isom.Hom.Ext}
Let $S$ be a (smooth) surface, embedded in $\mathbf P^N$, let $\mathcal J$ be the ideal sheaf of $S$ in $\mathbf P^N$  and consider the connecting homomorphism
\begin{equation*}
\mathrm{Hom}(\mathcal J/\mathcal J^2,\omega_{S}(-1)) \overset{\delta} \longrightarrow \mathrm{Ext}^1(\Omega_{S},\omega_{S}(-1)).
\end{equation*}
\begin{enumerate}
\item If $p_g(S)=0$ and $h^1(\mathcal O_{S}(1))=0$, then $\delta$ is injective;
\item if $q(S)=0$ and $S$ is embedded by a complete linear series,
then $\delta$ is surjective.
\end{enumerate}
\end{lemma}

\begin{proof}
The map $\delta$ fits in the following exact sequence:
\begin{equation*}
\mathrm{Hom}(\Omega_{\mathbf P^N}|_S,\omega_S(-1)) \longrightarrow \mathrm{Hom}(\mathcal J/\mathcal J^2,\omega_S(-1)) \overset{\delta} \longrightarrow \mathrm{Ext}^1(\Omega_S,\omega_S(-1))  \longrightarrow \mathrm{Ext}^1(\Omega_{\mathbf P^N}|_S,\omega_S(-1)).
\end{equation*}
Using the restriction to $S$ of the Euler sequence we get
\begin{equation*}
 \mathrm{Hom}(\mathcal O_S^{N+1}(-1),\omega_S(-1)) \longrightarrow \mathrm{Hom}(\Omega_{\mathbf P^N}|_S,\omega_S(-1)) \longrightarrow \mathrm{Ext}^1(\mathcal O_S,\omega_S(-1)).
 \end{equation*}
 Thus, if $p_g(S)=h^1(\mathcal O_S(1))=0$
 then $\mathrm{Hom}(\Omega_{\mathbf P^N}|_S,\omega_S(-1))$ vanishes so $\delta$ is injective. On the other hand, from the Euler sequence we also obtain
 \begin{equation*}
\mathrm{Ext}^1(\mathcal O_S^{N+1}(-1),\omega_S(-1)) \to \mathrm{Ext}^1(\Omega_{\mathbf P^N}|_S,\omega_S(-1)) \to \mathrm{Ext}^2(\mathcal O_S,\omega_S(-1)) \overset{\gamma} \to \mathrm{Ext}^2(\mathcal O_S^{N+1}(-1),\omega_S(-1)),
\end{equation*}
where $\gamma$ is the dual of
\begin{equation*}
H^0(\mathcal O_{\mathbf P^N}(1)) \otimes H^0(\mathcal O_S) \overset{\gamma^*} \longrightarrow H^0(\mathcal O_{S}(1)).
\end{equation*}
Then, if $S$ is embedded by a complete linear series, $\gamma^*$ is surjective and $\gamma$ is injective. If, in addition, $q(S)=0$, then $\mathrm{Ext}^1(\Omega_{\mathbf P^N}|_S,\omega_S(-1))$ vanishes and $\delta$ is surjective.
\end{proof}

\begin{theorem}\label{exist.embedd.carpets}
With the notation of~\ref{setup} and~\ref{notation.blowups}, let
 $i$ be the embedding of $Y$ in $\mathbf P^N$ induced by $|dL -E|$, where $(d,s)$ satisfies condition (a), (b), (c) or (d) of Proposition~\ref{blowup.embeddings}.
\begin{enumerate}
 \item[(1.1)] If $d=2$ and $s=1$, or
\item[(1.2)] if $d \geq 3$ and $s \leq \frac{1}{2}d^2-\frac{1}{2}d+1$,
\end{enumerate}
then $\mathrm{Hom}(\mathcal I/\mathcal I^2, \omega_Y(-1)) = 0$; thus, there do not exist double structures on $i(Y)$, embedded in $\mathbf P^N$, whose conormal line bundle is a subsheaf of $\omega_Y(-1)$.
\smallskip

\noindent On the contrary,
\begin{enumerate}
\item[(2.1)] if $d=3$ and $s=5,6$, or
\item[(2.2)] if $d=4$ and $s=8,9,10$, or
\item[(2.3)] if $d \geq 5$ and $\frac{1}{2}d^2-\frac{1}{2} < s \leq \frac{1}{2}d^2+\frac{3}{2}d-5$,
\end{enumerate}
then $\mathrm{Hom}(\mathcal I/\mathcal I^2, \omega_Y(-1)) \neq 0$; thus, there exist double structures on $i(Y)$, embedded in $\mathbf P^N$, whose conormal line bundle is a subsheaf of $\omega_Y(-1)$. In fact, every double structure $\tilde Y$ on $Y$ with conormal bundle $\omega_Y(-1)$ admits a unique morphism $\tilde i$ to $\mathbf P^N$ extending $i$. The image of $\tilde i$ is a double structure on $i(Y)$ whose conormal line bundle is a subsheaf of $\omega_Y(-1)$, unless $\tilde Y$ is the split double structure, in which case $\tilde i(\tilde Y)=i(Y)$.

\end{theorem}

\begin{proof}
First note that, since $x_1,\dots,x_s$ are in general position, if $(d,s)$ satisfies (a), (b), (c) or (d) of Proposition~\ref{blowup.embeddings}, then
\begin{equation*}
H^1(\mathcal O_{\mathbf P^2}(d) \otimes \mathfrak m)=0,
\end{equation*}
 so, by the  projection formula and the Leray spectral sequence,
\begin{equation}\label{vanish.for.isom.Hom.Ext}
 H^1(\mathcal O_Y(1))=0.
\end{equation}
On the other hand, $p_g(Y)=q(Y)=0$ and $i$ is induced by a complete linear series,
 so we can apply Lemma~\ref{isom.Hom.Ext} and deduce that
the connecting homomorphism
\begin{equation*}
\mathrm{Hom}(\mathcal I/\mathcal I^2,\omega_{Y}(-1)) \overset{\delta} \longrightarrow \mathrm{Ext}^1(\Omega_{Y},\omega_{Y}(-1))
\end{equation*}
is an isomorphism. Then the result follows from Corollary~\ref{exist.abstr.carpets} and~\cite[Proposition 2.1, (2)]{Gon} (see also~\cite[Lemma 1.4]{GPcarpets} or~\cite{HV}).
\end{proof}

\noindent Before going on our way to study the deformations of the canonical map of the canonical double covers of $Y$ constructed in Corollary~\ref{exist.emb.cover.cor} we will need to introduce some extra notation and conventions:

\begin{noname}\label{setup.extra}
{\bf Set--up and notation.}
 {\rm With addition to the notations introduced in~\ref{setup} and~\ref{notation.blowups}, throughout  the remaining of the article we will assume
\begin{enumerate}
 \item $s$ and $d$ to be as in (a), (b), (c), (d) or (e) of Proposition~\ref{exist.emb.cover};
\item $B$ to be a smooth divisor in $|-2K_Y+2dL-2E|$ (see Proposition~\ref{exist.emb.cover});
\item $\pi:X \longrightarrow Y$ to be the double cover of $Y$ branched along $B$;
\item $\varphi$ to be, as in Corollary~\ref{exist.emb.cover.cor}, the canonical map of $X$;
\item $i$ to be, as in Theorem~\ref{exist.embedd.carpets}, the embedding of $Y$ in $\mathbf P^N$ induced by $|dL -E|$ (see Proposition~\ref{exist.emb.cover}).
\end{enumerate}}
\end{noname}

\begin{remark}\label{remark.canon.morphism}
{\rm Recall that, in the previous set--up~\ref{setup.extra}, by Corollary~\ref{exist.emb.cover.cor}, $X$ is a surface of general type and $\varphi$ is a morphism to $\mathbf P^N$ and is the composition of $\pi$ followed by $i$.}
\end{remark}

\noindent In order to complete the knowledge of $\Psi_2$ needed to apply Theorem~\ref{Psi2<>0} and also to obtain the conditions about unobstructedness required in this theorem, we look now at the morphisms $\pi$ and $\varphi$. Recall that  according to~\ref{setup.extra}, $s$ and $d$ are as in (a), (b), (c), (d) or (e) of Proposition~\ref{exist.emb.cover}. On the other hand, by Theorem~\ref{exist.embedd.carpets}, $\Psi_2$ could be non zero only if $s > \frac{1}{2}d^2-\frac{1}{2}$, i.e, only if $(d,s)=(3,5), (3,6), (4,8), (4,9), (4,10), (5,13)$ or $(5,14)$. These values of $(d,s)$ are dealt with in the next

\begin{proposition}\label{vanish.normal} With the notation of~\ref{setup},~\ref{notation.blowups} and~\ref{setup.extra},
if $(d,s)=(3,5), (3,6), (4,8), (4,9), (4,10),$ $(5,13)$ or $(5,14)$,
then
\begin{enumerate}
\item $H^1(\mathcal N_\pi)=0$; and
\item $H^1(\mathcal N_\varphi)=0$.
\end{enumerate}
\end{proposition}

\begin{proof}
\noindent{Proof of (1)}. By~\eqref{cohom.normalpi.formula} $H^1(\mathcal N_\pi) \simeq H^1(\mathcal O_B(B))$. To compute the latter consider
\begin{equation}\label{proving.h1normal.zero}
H^1(\mathcal O_Y) \longrightarrow H^1(\mathcal O_Y(B)) \longrightarrow H^1(\mathcal O_B(B)) \longrightarrow H^2(\mathcal O_Y).
\end{equation}
By Serre duality, $h^2(\mathcal O_Y)=p_g(Y)=0$ and recall also that $h^1(\mathcal O_Y)=q(Y)=0$, so $H^1(\mathcal O_B(B))=H^1(\mathcal O_Y(B))$. By~\eqref{canon.hyper2}, $\mathcal O_Y(B)=p^*\mathcal O_{\mathbf P^2}(2d+6) \otimes \mathcal O_Y(-4E)$, so by the projection formula and the Leray spectral sequence,
$H^1(\mathcal O_Y(B))=H^1(\mathcal O_{\mathbf P^2}(2d+6) \otimes \mathfrak m^{4})$. Now we see that
$H^1(\mathcal O_{\mathbf P^2}(2d+6) \otimes \mathfrak m^{4})$ is the cokernel of the evaluation map
\begin{equation}\label{beta.restriction.to.4S}
H^0(\mathcal O_{\mathbf P^2}(2d+6)) \overset{\beta} \longrightarrow H^0(\mathcal O_{\mathbf P^2}(2d+6)|_\Sigma),
\end{equation}
where $\Sigma$ is the third infinitesimal neighborhood of $x_1,\dots,x_s$.
Now we want to apply~\cite[Theorems 2.4 and 5.2]{CilibertoMiranda}. Note first that if $d \geq 5$, then $s > \frac{1}{2}d^2-\frac{1}{2} \geq 12$, so~\cite[Theorem 2.4]{CilibertoMiranda} implies that $|(2d+6)L-4E|$ is not $(-1)$--special. Then~\cite[Theorem 5.2]{CilibertoMiranda} and the fact that
$|(2d+6)L-4E| \neq \emptyset$ imply the vanishing of $H^1(\mathcal O_{\mathbf P^2}(2d+6) \otimes \mathfrak m^{4})$ in this case.
If $d=4$, then $s > \frac{1}{2}d^2-\frac{1}{2}$ implies $s \geq 8$. Then a simple computation
shows that~\cite[Theorems 2.4 and 5.2]{CilibertoMiranda} and the fact that
$|(2d+6)L-4E| \neq \emptyset$ imply as before the vanishing of $H^1(\mathcal O_{\mathbf P^2}(2d+6) \otimes \mathfrak m^{4})$.
Finally~\cite[Theorem 2.4]{CilibertoMiranda} shows that there do not exist homogeneous, $(-1)$--special systems of degree $3$ and multiplicity $4$, so if $d=3$, ~\cite[Theorem 5.2]{CilibertoMiranda} and the fact that
$|(2d+6)L-4E| \neq \emptyset$ imply the vanishing of $H^1(\mathcal O_{\mathbf P^2}(2d+6) \otimes \mathfrak m^{4})$. Thus $H^1(\mathcal O_B(B))$, and therefore, $H^1(\mathcal N_\pi)$ equals $0$.

\medskip

\noindent{Proof of (2)}. Recall the sequence (see~\cite[(3.3.2)]{Gon})
\begin{equation*}
0 \longrightarrow \mathcal N_\pi \longrightarrow \mathcal N_\varphi \longrightarrow \pi^*\mathcal N_{i(Y),\mathbf P^N} \longrightarrow 0.
\end{equation*}
Since $ \mathcal N_\pi$ is supported on the curve $R$ we have the short exact sequence
\begin{equation}\label{isom2}
 H^1(\mathcal N_\pi) \longrightarrow H^1(\mathcal N_\varphi) \longrightarrow H^1(\mathcal N_{Y,\mathbf P^N}) \oplus H^1(\mathcal N_{Y,\mathbf P^N} \otimes \omega_Y(-1)) \longrightarrow 0
\end{equation}

\smallskip

\noindent We prove now that both $H^1(\mathcal N_{Y,\mathbf P^N})$ and $H^1(\mathcal N_{Y,\mathbf P^N} \otimes \omega_Y(-1))$ vanish. The normal sequence and the Euler sequence of $\mathbf P^N$ say that, in order to prove the vanishing of $H^1(\mathcal N_{Y,\mathbf P^N})$, it will suffice to check the vanishings of $H^1(\mathcal O_Y(1))$, $H^2(\mathcal O_Y)$ and $H^2(\mathcal T_Y)$. Now recall that $h^2(\mathcal O_Y)=p_g(Y)=0$. On the other hand, $H^1(\mathcal O_Y(1))=0$ (see~\eqref{vanish.for.isom.Hom.Ext}). We deal finally with $H^2(\mathcal T_Y)$. Recall that
$H^2(\mathcal T_Y)=H^0(\Omega_Y \otimes \omega_Y)^{\vee}$. Now to compute $H^0(\Omega_Y \otimes \omega_Y)$
we use~\eqref{rel.cotangent.seq} and get
\begin{equation*}
0 \longrightarrow H^0(p^*\Omega_{\mathbf P^2}\otimes \omega_Y) \longrightarrow H^0(\Omega_Y \otimes \omega_Y) \longrightarrow   H^0(\Omega_{Y/\mathbf P^2} \otimes \omega_Y).
\end{equation*}
To prove the vanishing of $H^0(p^*\Omega_{\mathbf P^2} \otimes \omega_Y)$ we use that, by the Euler sequence, it is contained in
$H^0(p^*\mathcal O_{\mathbf P^2}(-4) \otimes \mathcal O_Y(E))^{\oplus 3}$, which is $0$.
On the other hand, by~\eqref{rel.cotangent}, $\Omega_{Y/\mathbf P^2} \otimes \omega_Y=\mathcal O_{E_1}(3E_1) \oplus \cdots \oplus \mathcal O_{E_s}(3E_s)$, which has no global sections. This implies that  \begin{equation}\label{h2.tangent.Y=0}
H^2(\mathcal T_Y)=0.
\end{equation}

\smallskip

\noindent Now we deal with $H^1(\mathcal N_{Y,\mathbf P^N} \otimes \omega_Y(-1))$. We consider the sequence
\begin{equation}\label{rel.cotangent.seq2}
H^1(\mathcal T_{\mathbf P^N|_Y} \otimes \omega_Y(-1)) \longrightarrow  H^1(\mathcal N_{Y,\mathbf P^N} \otimes \omega_Y(-1)) \longrightarrow H^2(\mathcal T_Y\otimes \omega_Y(-1)).
\end{equation}
By Serre's duality $H^2(\mathcal T_Y\otimes \omega_Y(-1))=H^0(\Omega_Y(1))^{\vee}$. Now we twist and take cohomology on
\eqref{rel.cotangent.seq} to prove  $H^0(\Omega_Y(1))=0$. We have $H^0(\mathcal O_{E_i}(E_i))=H^0(\mathcal O_{\mathbf P^1}(-1))=0$, so we just need to show that $H^0(p^*\Omega_{\mathbf P^2} \otimes \mathcal O_Y(1))=0$. This group can be identified with the kernel of the multiplication map $\alpha$ defined in Lemma~\ref{surjectivity}. We check the injectivity of $\alpha$ case by case. If $(d,s)=(3,6)$ or $(4,10)$,
then $H^0(\mathcal O_{\mathbf P^2}(d-1) \otimes \mathfrak m)=0$, so $\alpha$ is injective. If $(d,s)=(3,5),  (4,9)$ or $(5,14)$,
then $h^0(\mathcal O_{\mathbf P^2}(d-1) \otimes \mathfrak m)=1$, so $\alpha$ is also injective. Finally, if $(d,s)=(4,8)$ or $(5,13)$, then $h^0(\mathcal O_{\mathbf P^2}(d-1) \otimes \mathfrak m)=2$.
In this case we look at the proof of Lemma~\ref{surjectivity} and at diagram~\eqref{restrict.tocurves}. From that we see that $\alpha$ is injective if $\alpha'$ is injective. In this case ($(d,s)=(4,8)$ or $(5,13)$), $\alpha_1'$ is injective, so $\alpha'$ is injective if $\alpha_2'$.
Now the injectivity of $\alpha_2'$ follows from the base--point--free pencil trick (see~\cite[p. 126]{ACGH}).
Thus
\begin{equation}\label{tangent.vanish}
H^2(\mathcal T_Y\otimes \omega_Y(-1))=0
\end{equation}
To study $H^1(\mathcal T_{\mathbf P^N|_Y} \otimes \omega_Y(-1))$ we use again the Euler sequence of $\mathbf P^N$. After restricting the Euler sequence to $Y$, twisting  and taking cohomology, we get
\begin{equation*}
H^1(\omega_Y)^{\oplus N+1} \longrightarrow H^1(\mathcal T_{\mathbf P^N|_Y} \otimes \omega_Y(-1))  \longrightarrow H^2(\omega_Y(-1)) \overset{\rho^*}  \longrightarrow H^2(\omega_Y)^{\oplus N+1}.
\end{equation*}
Now, $H^1(\omega_Y)=H^1(\mathcal O_Y)^{\vee}=0$, for $q(Y)=0$. On the other hand, $\rho^*$ is dual of the multiplication map of global sections
\begin{equation*}
H^0(\mathcal O_Y(1)) \otimes H^0(\mathcal O_Y) \overset{\rho} \longrightarrow H^0(\mathcal O_Y(1)),
\end{equation*}
which is obviously an isomorphism. Thus $H^1(\mathcal T_{\mathbf P^N|_Y} \otimes \omega_Y(-1))=0$.  This together with~\eqref{rel.cotangent.seq2} and~\eqref{tangent.vanish} gives the vanishing of $H^1(\mathcal N_{Y,\mathbf P^N} \otimes \omega_Y(-1))$.  This vanishing and the vanishing of $H^1(\mathcal N_{Y,\mathbf P^N})$,  having in account~\eqref{isom2} and (1), complete the proof of (2).
\end{proof}

\noindent Finally, we have all the tools needed for using Theorems~\ref{Psi2=0} and~\ref{Psi2<>0} to see how the deformations of the double canonical covers of $Y$ are. We deal with the two possible behaviors (either $\varphi$ can be deformed to a degree $1$ morphism or $\varphi$ deforms always to a degree $2$ morphism) in two separate theorems, Theorems~\ref{theorem.construct.surfaces} and~\ref{2to1deformsto2to1}:

\begin{theorem}\label{theorem.construct.surfaces} Let $(d,s)=(3,5), (3,6), (4,8), (4,9), (4,10), (5,13)$ or $(5,14)$.
Then $X$ is unobstructed. Moreover,  there exist a smooth  irreducible algebraic curve $T$ and a point $0 \in T$ such that the morphism $\varphi: X \longrightarrow \mathbf P^N$ (which is finite and $2:1$ onto $i(Y)$) fits into a flat family $\Phi: \mathcal X \longrightarrow \mathbf P^N_T$ over $T$, satisfying
\begin{enumerate}
\item $\Phi_0=\varphi$;
\item For any $t \neq 0$, $t \in T$, $\mathcal X_t$ is a smooth and irreducible surface of general type, $\Phi_t: \mathcal X_t \longrightarrow \mathbf P^N_t$ is the canonical map of $\mathcal X_t$ and is finite of degree $1$ onto its image.
\end{enumerate}
Thus the canonical map of a surface corresponding to a general point of the component of $X$ in its moduli space is a finite morphism of degree $1$.
\end{theorem}

\begin{proof}
Recall (see Remark~\ref{remark.canon.morphism}) that $X$ is a surface of general type and $\varphi$ is its canonical map, which is a finite morphism of degree $2$. The surface $X$ is unobstructed by Lemma \ref{varphi.iff.X.unobstructed}, Remark \ref{canonical.remark}, (3) and Proposition~\ref{vanish.normal}, (2).
Then the result follows from Theorem \ref{Psi2<>0}, Theorem~\ref{exist.embedd.carpets}, (2) and Proposition~\ref{vanish.normal}, (1).
\end{proof}

\noindent The interest of the families of examples constructed in Theorem~\ref{theorem.construct.surfaces} lies in the fact that, 
apart from the case $(d,s)=(3,5)$  and maybe $(4,8)$ (see Question~\ref{moduli4,8}),
they are not complete intersections. Moreover, as we see in Proposition~\ref{noin3folds}, they cannot be obtained as smooth divisors of smooth rational scrolls of dimension $3$ (an easy way of producing canonically embedded surfaces, by adjunction). Thus the cases $(d,s)=(4,9), (4,10), (5,13)$ and $(5,14)$ (the case $(d,s)=(3,6)$ appears in~\cite{Horikawa.small3}) provide new, interesting examples of surfaces with finite birational canonical morphisms to low dimensional projective spaces which cannot be constructed by simple methods (in fact, surfaces with the same invariants as them cannot be either complete intersections or divisors in smooth rational scrolls of dimension $3$, as we will see in Proposition~\ref{noin3folds}). Especially interesting is the case $(d,s)=(4,10)$, for it exhibits a family of surfaces (necessarily with non very ample canonical line bundle) whose canonical map  is a $1:1$ morphism to $\mathbf P^4$. We make all these points precise in 
Remark~\ref{remarkd3} and Proposition~\ref{noin3folds}:

\begin{remark}\label{remarkd3} Assuming the hypothesis of Theorem~\ref{theorem.construct.surfaces},
\begin{enumerate}
 \item if $d=3$ and $s=5$, then the canonical map of a surface corresponding to a general point of the component of $X$ in its moduli space is an embedding into $\mathbf P^4$;
\item if $d=3$ and $s=6$, then the canonical map of a surface $X'$ corresponding to a  general point of the component of $X$ in its moduli space is not an embedding but maps $X'$ onto a sextic surface in $\mathbf P^3$, singular along a plane cubic curve;
\item if $d=4$ and $s=10$, then the canonical map of a surface $X'$ corresponding to a  general point of the component of $X$ in its moduli space is not an embedding but maps $X'$ onto a singular surface in $\mathbf P^4$.
\end{enumerate}
\end{remark}

\begin{proof}
If $d=3$ and $s=5$, then $p_g=5$ and $Y$ is a smooth rational surface of degree $4$ in $\mathbf P^4$, which is therefore the complete intersection of two quadrics $Q_1$ and $Q_2$. Then $\mathcal I /\mathcal I^2 \simeq \mathcal O_Y(-2) \oplus \mathcal O_Y(-2)$ and $\omega_Y(-1) \simeq \mathcal O_Y(-2)$, so obviously there are surjective homomorphisms in $\mathrm{Hom} (\mathcal I/\mathcal I^2, \omega_Y(-1))$ (recall that $\mathrm{Hom}(\mathcal I/\mathcal I^2, \omega_Y(-1))$ contains a surjective homomorphism if and only if there exists a double structure on $Y$, embedded in $\mathbf P^N$, with conormal bundle $\omega_Y(-1)$; then, in our case one of such surjective homomorphisms would correspond to the double structure on $Y$ obtained by considering the intersection of $Q_1$ and the (unique) double structure on $Q_2$ inside $\mathbf P^4$).
Then~\cite[Theorem 1.4]{infinitesimal} implies that $\varphi$ can be deformed to an embedding. Obviously, a smooth complete intersection of a quadric and a quartic threefold is the image of one of such canonical embeddings.

\smallskip
\noindent If $d=3$ and $s=6$, then $p_g=4$ and $c_1^2=6$.  If the canonical morphism of $X'$ were an embedding, $X'$ would be a canonically embedded smooth surface in $\mathbf P^3$, but it would have degree $6$, which would be impossible. The last claim of the statement was proved by Horikawa  (see~\cite[Theorem 3.2]{Horikawa.small3}).

\noindent If $d=4$ and $s=10$ and the canonical morphism of $X'$ were an embedding, then there would exist a smooth surface in $\mathbf P^4$, of degree $12$ and canonically embedded. This is impossible by the double point formula (see Example 4.1.3 of the Appendix A of~\cite{Hart}). The impossibility for $X'$ to be canonically embedded is also suggested by the fact that in this case  $\mathrm{Hom} (\mathcal I/\mathcal I^2, \omega_Y(-1))$ does not contain surjective homomorphisms, because $c_2(\mathcal N_{Y/\mathbf P^4} \otimes \omega_Y(-1)) \neq 0$.

\end{proof}

\noindent Before stating the next proposition, we recall the following notation:

\begin{noname}\label{scroll.notation}
{\bf Notation.}
{\rm Let $a, b$ and $c$ be natural numbers such that $a \leq b \leq c$ and let $r=a+b+c+3$. Then $S(a,b,c)$ is the smooth rational normal scroll $Z$ of dimension $3$ and degree $r-3$ in $\mathbf P^{r-1}$ obtained as the union of the spans of triplets of corresponding points lying on three rational normal curves of degrees $a, b$ and $c$ in $\mathbf P^a,  \mathbf P^b$ and $\mathbf P^c$ respectively. We call $H$ the hyperplane divisor of $Z$ and  $F$ the fiber of the projection of $Z$ to $\mathbf P^1$.}
\end{noname}

\medskip

\begin{proposition}~\label{noin3folds}
Let $m$ and $l$ be integers with $m \geq 4$ and such that $mH+lF$ is a base--point--free and big divisor in $Z=S(a,b,c)$.
Let $S$ be a smooth divisor of the  linear system $|mH+lF|$. The pair $(p_g(S),c_1^2(S))$ cannot be the pair of invariants (see Proposition~\ref{cancover.blownP2.invariants}) of any surface $X$ in Theorem~\ref{theorem.construct.surfaces} or of any surface with finite birational canonical morphism as the ones constructed in Theorem~\ref{theorem.construct.surfaces}, except if $(p_g(S),c_1^2(S))=(5,8)$; in this case there exist surfaces $S$ as above with very ample canonical divisor. In particular the surfaces with finite birational canonical morphism constructed in Theorem~\ref{theorem.construct.surfaces}, except maybe the ones with $(p_g,c_1^2)=(5,8)$, are not smooth divisors in smooth rational normal scrolls of dimension $3$.
\end{proposition}

\begin{proof}
Suppose first that $S$ is a smooth surface in $|mH+lF|$. Then $K_S=(K_Z+S)|_S$ and $c_1^2(S)=(K_Z+S)^2 \cdot S$. Since $h^0(K_Z)=h^1(K_Z)=0$, then $p_g(S)=h^0(K_Z+S)$. Recall that $K_Z=-3H+(r-5)F$; then we can compute $h^0(K_Z+S)$ using the Hirzebruch--Riemann--Roch formula and we can write $p_g(S)$ and $c_1^2(S)$ in terms of $r, m$ and $l$:
\begin{align}\label{parametric}
p_g(S)&=\frac 1 6 (m-2)(m-1)(rm+3l)-\frac 1 2 (m-2)(m-1)(m+1) \cr
c_1^2(S)&=(m-3)(m-1)(rm+3l)-m(m-3)(3m+1).
\end{align}
It is possible to eliminate $r$ and $l$ from the above equations and conclude that the pair $(p_g(S),c_1^2(S))$ is a point $(x',y) \in \mathbf N^2$ lying on the line of equation
\begin{equation}\label{implicit}
y=6\frac{m-3}{m-2}x'-(m-3)(m+3),
\end{equation}
where $m$ is an integer, $m \geq 4$. Consider now the pairs of invariants $(p_g,c_1^2)$ corresponding to the surfaces $X$ in Theorem~\ref{theorem.construct.surfaces} or to surfaces with finite birational canonical morphism constructed in Theorem~\ref{theorem.construct.surfaces}. They are $(4,6), (5,8), (5,12), (6,14), (7,16)$, $(7,22)$ and $(8,24)$ (see Proposition~\ref{cancover.blownP2.invariants}). A routine computation yields that the only one among those pairs lying on one of the lines~\eqref{implicit} is the pair $(5,8)$ if $m=4$.
In this case it is possible to find smooth surfaces $S$ in $|4H+lF|$
just by setting $l=-4$, $a=1$ and $b=c=2$. In addition, $S$ can be chosen so that $K_S$ is very ample. Indeed, $K_Z+S=H-F$ is base--point--free (although not ample).  Let $C$ be the only curve in $Z$ contracted by $H-F$.
Then $K_Z+S$ is very ample on $Z \smallsetminus C$. On the other hand, since $4H+lF$ is also base--point--free and $(4H-4F)\cdot C=0$, we can choose $S \subset Z \smallsetminus C$, so $K_S=(K_Z+S)|_S$ is very ample.
\end{proof}

\noindent Finally we find out the cases for which the canonical double cover of $Y$ only deforms to a degree $2$ morphism:

\begin{theorem}\label{2to1deformsto2to1} Let $d$ and $s$ be such that
\begin{enumerate}
\item $d=2$ and $s=1$; or
 \item $3 \leq d \leq 6$ and $s \leq \frac{1}{2}d^2-\frac{1}{2}d+1$; or
\item $d \geq 7$ and $s \leq \frac{1}{5}d^2+\frac{13}{10}d+\frac{21}{10}$.
\end{enumerate}
Then $X$ is unobstructed. Moreover, for any deformation $\Phi: \mathcal X \longrightarrow \mathbf P^N_T$ of $\varphi: X \longrightarrow \mathbf P^N$ over a smooth irreducible algebraic curve $T$, $\mathcal X_t$ is a regular surface of general type and $\Phi_t: \mathcal X_t \longrightarrow \mathbf P^N_t$ is its canonical map, which is a finite morphism of degree $2$.
Thus the canonical map of a surface corresponding to a general point of the component of $X$ in its moduli space is a finite morphism of degree $2$.
\end{theorem}

\begin{proof}
We want to apply Theorem~\ref{Psi2=0}.
Corollary~\ref{exist.emb.cover.cor} tells us that $X$ is a smooth surface of general type whose canonical map $\varphi$ is a finite morphism of degree $2$. Recall that $q(Y)=0$. We have also $H^1(\mathcal O_Y(1))=0$ (see ~\eqref{vanish.for.isom.Hom.Ext}) and, since $p_g(Y)=0$, $h^0(\omega_Y(-1))=0$. We check also that $H^1(\omega_Y^{-2}(2))=0$. For this  we use~\cite[Theorems 2.4 and 5.2]{CilibertoMiranda}. Recall that $H^1(\omega_Y^{-2}(2))=H^1(\mathcal O_{\mathbf P^2}(2d+6) \otimes \mathfrak m^4)$, by the projection formula and the Leray spectral sequence.
Going through the statement of~\cite[Theorem 2.4]{CilibertoMiranda} one realizes that $|(2d+6)H-4E|$ would be $(-1)$--special only if $s=2$ and $2d+6=4,5$ or $6$; or if $s=3$ and $2d+6=6$; or if $s=5$ and $2d+6=8$ or $9$. In the three cases, $d \leq 1$, so $|(2d+6)H-4E|$ is not $(-1)$--special. Then~\cite[Theorem 5.2]{CilibertoMiranda} and the fact that
$|(2d+6)H-4E| \neq \emptyset$ implies the vanishing of $H^1(\mathcal O_{\mathbf P^2}(2d+6) \otimes \mathfrak m^{4})$.
On the other hand, note that the arguments used in the proof of Proposition~\ref{vanish.normal} to show that $H^2(\mathcal T_Y)=0$, work also under the present hypothesis, so $Y$ is unobstructed. Finally Theorem~\ref{exist.embedd.carpets}, (1) implies $\Psi_2=0$.
Then we can apply Theorem~\ref{Psi2=0}.
and the result follows.
\end{proof}

\begin{question}\label{question}
{\rm If $d=5$ and $s=12$
or if $d=6$ and $s=17$, then, according to Corollary~\ref{exist.emb.cover.cor}, there exist smooth surfaces of general type whose canonical map $\varphi$ is a morphism of degree $2$ onto $Y$, embedded by $|dH-E|$. For these surfaces the question of whether $\varphi$ can be deformed to a morphism of degree $1$ remains open.}
\end{question}

\section{Consequences for the topology, geography and moduli}\label{moduli.section}

In this section we compute the invariants of the surfaces we have constructed in Section~\ref{construct.from.blownup.section}. This way we pinpoint their position in  the geography of surfaces of general type. In addition, we compute the dimension of the moduli components parameterizing our surfaces, as well as the codimension of the loci parameterizing surfaces whose canonical map is a degree $2$ morphism onto its image. We also discover interesting components of the moduli space of surfaces of general type: in Theorem~\ref{2components} we show that, for infinitely many moduli spaces of surfaces of general type, there exist at least two components; the general point of one of the two corresponds to a surface whose canonical map is a degree $2$ morphism whereas the general point of the other corresponds to a surface that can be canonically embedded. We use the same notations and conventions given in~\ref{setup}, \ref{notation.blowups} and~\ref{setup.extra}.

\begin{proposition}\label{cancover.blownP2.invariants}
1) The surfaces of general type appearing in Theorem~\ref{theorem.construct.surfaces} have the following invariants:

\bigskip

\centerline{\vbox{\tabskip=0pt \offinterlineskip
%\tabskip= .25 truecm
\def\tablerule{\noalign{\hrule}}
\halign to 8truecm
%{\valign to125pt}
{\strut
#& \vrule#
\tabskip=0em plus 3em
%\tabskip= .25 truecm
&
\hskip .3cm
\hfil# \hskip .1cm
& \vrule #
& \hskip .2cm
\hfil #
\hfil  \hskip .1cm
& \vrule \vrule \vrule#&
\hskip .2cm \hfil#\hfil \hskip .1cm & \vrule#&
\hskip .2cm \hfil#\hfil \hskip .1cm & \vrule#&
\hskip .2cm \hfil#\hfil \hskip .1cm & \vrule#&
\hskip .2cm \hfil#\hfil \hskip .1cm &
\vrule#&
\hskip .25cm \hfil# %\hskip .05truecm
\hfil & \vrule#\tabskip=0pt
\cr\tablerule
%&&Type
&&%\omit\hidewidth
$d$ %\hidewidth
&& %\omit\hidewidth
$s$ %\hidewidth
&&%\omit\hidewidth
$p_g$ %\hidewidth
&&%\omit\hidewidth
$q$ %\hidewidth
&&%\omit\hidewidth
 $\chi$  %\hidewidth
&&%\omit\hidewidth
  $c_1^2$  %hidewidth
&&%\omit\hidewidth
$c_1^2/c_2$ %\hidewidth
&\cr\tablerule
%\cr
\tablerule
\tablerule
&&$3$
&&$6$
&&$4$
&&$0$
&&$5$
&&$6$
&&$1/9$
&\cr
%&\vskip -2truecm&\vskip -2truecm&\vskip -2truecm&\vskip -2truecm&\vskip -2truecm&\vskip -2truecm&\vskip -2truecm&\vskip -2truecm&\vskip -2truecm&\vskip -2truecm&\vskip -2truecm&\vskip -2truecm&\vskip -2truecm&\vskip -2truecm&\cr
\tablerule
%&&2
&&$3$
&&$5$
&&$5$
&&$0$
&&$6$
&&$8$
&&$1/8$
&\cr
\tablerule
%&&3
&&$4$ %\ *$ \hskip -.2truecm
&&$10$
&&$5$
&&$0$
&&$6$
&&$12$
&&$1/5$
&\cr
\tablerule
%&&4
&&$4$ % \ *$ \hskip -.2truecm
&&$9$
&&$6$
&&$0$
&&$7$
&&$14$
&&$1/5$
&\cr
\tablerule
%&&5.1
&&$4$ %\hskip -.15truecm
&&$8$
&&$7$
&&$0$
&&$8$
&&$16$
&&$1/5$
&\cr
\tablerule
&&$5$ %\hskip -.15truecm
&&$14$
&&$7$
&&$0$
&&$8$
&&$22$
&&$11/37$
&\cr
\tablerule
%&&6.1
&&$5$ %\hskip -.15truecm
&&$13$
&&$8$
&&$0$
&&$9$
&&$24$
&&$2/7$
&\cr
\tablerule
\noalign{\smallskip}}}}
\medskip

2) The surfaces of general type appearing in Theorem~\ref{2to1deformsto2to1}
have the following invariants:
\begin{eqnarray*}
 p_g&=&\frac{1}{2}d^2+\frac{3}{2}d-s+1 \cr
 \cr
q&=&0 \cr
\cr
\chi&=&\frac{1}{2}d^2+\frac{3}{2}d-s+2 \cr
\cr
c_1^2&=&2d^2-2s \cr
\cr
\frac{c_1^2}{c_2}&=&\frac{d^2-s}{2d^2+9d-5s+12}
\end{eqnarray*}
\end{proposition}

\begin{proof}
The values for $p_g$ follow from the fact that, since $\{x_1,\dots,x_s\}$ are in general position, they impose independent conditions on curves of $\mathbf P^2$ of degree $d$. By the construction of $X$ (see (\ref{setup.extra})), $h^1(\mathcal O_X)=h^1(\mathcal O_Y) + h^1(\omega_Y(-1))$. Recall that $q(Y)=0$. On the other hand, since  $(d,s)$ satisfy (a), (b), (c) or (d) of Proposition~\ref{blowup.embeddings}, we have $h^1(\mathcal O_Y(1))=0$ (see~\eqref{vanish.for.isom.Hom.Ext}).
Thus $q(X)=0$.
The values of $c_1^2$ follow from the values of the degree of $Y$ inside $\mathbf P^N$, having in account that $\omega_X=\varphi^*\mathcal O_Y(1)$ and $\varphi$ has degree $2$ onto $Y$. Finally, the values of $\frac{c_1^2}{c_2}$ follow from Noether's formula.
\end{proof}

\begin{remark}\label{cancover.blownP2.geography}
{\rm We now present more graphically the information given in Proposition~\ref{cancover.blownP2.invariants}, by displaying on a plane the pairs $(x,y)=(\chi,c_1^2)$ of the surfaces of general type appearing in Theorems~\ref{theorem.construct.surfaces} and~\ref{2to1deformsto2to1}. Let us fix an integer $d$, $d \geq 2$. Then the points $(x,y)$ corresponding to the surfaces of Theorems~\ref{theorem.construct.surfaces} and~\ref{2to1deformsto2to1} are points (with integer coordinates) lying on the line $y-2x-(d^2-3d-4)=0$. Observe also that this line is parallel to and above the Noether's line, which is $y-2x+6=0$, except in the case $(d,s)=(2,1)$, since in this case $Y$ is embedded as a surface of minimal degree, so, under these circumstances, $(x,y)$ belongs in the Noether's line.
Then in each line $y-2x-(d^2-3d-4)=0$, $d \in \mathbf Z$, $d \geq 2$, the surfaces of Theorems~\ref{theorem.construct.surfaces} and~\ref{2to1deformsto2to1} yield a finite number of integer points. Precisely these are the integer points $(x,y)$
\begin{enumerate}
\item[] with $x=6$ if $(x,y)$ is on the line $y-2x+6=0$,
\item[] with $5 \leq x \leq 10$ if $(x,y)$ is  on the line $y-2x+4=0$,
\item[] with $6 \leq x \leq 15$ if $(x,y)$ is  on the line $y-2x=0$,
\item[] with $8 \leq x \leq 21$ if $(x,y)$ is  on the line $y-2x-6=0$,
\item[] with $13 \leq x \leq 28$ if $(x,y)$ is  on the line $y-2x-14=0$,
\end{enumerate}
and if $d \geq 7$, these are the integer points $(x,y)$ that lie on the lines $y-2x-(d^2-3d-4)=0$ ($d \in \mathbf Z$) and in the region between the parabolas of equations $256x^2-96xy+9y^2-638x+44y=0$ and $16x^2-8xy + y^2 -48x-6y=0$.
}
\end{remark}

\noindent Next we remark that the surfaces constructed in Theorems~\ref{theorem.construct.surfaces} and~\ref{2to1deformsto2to1} are not only regular, but also simply connected:

\begin{remark}\label{simply.connected}
The surfaces of general type that we have constructed in Theorems~\ref{theorem.construct.surfaces} and~\ref{2to1deformsto2to1} are simply connected.
\end{remark}

 \begin{proof} Recall that $Y$ is $\mathbf P^2$  blown up at a finite number of points, hence the fundamental group of $Y$ is the same as the fundamental group of $\mathbf P^2$, which is $0$. The morphism $\pi$ is a double cover of $Y$ branched along a nef and big divisor (see Lemma~\ref{Harbourne} and the proof of Proposition~\ref{exist.emb.cover.cor}, where we prove that the branch divisor of our covers is base--point--free). Then  the fundamental group of $X$ is the same as the fundamental group of $Y$ by~\cite[Cor. 2.7]{Nori} (note that the ampleness hypothesis required there can be relaxed to big and nefness), so $X$ is simply connected.  Then,  consider the families of surfaces associated to the deformations of $X$ given in Theorems~\ref{theorem.construct.surfaces} and~\ref{2to1deformsto2to1}. All the smooth fibers in such families are diffeomorphic to each other, hence they are also simply connected. Thus the surfaces constructed in  Theorems~\ref{theorem.construct.surfaces} and~\ref{2to1deformsto2to1} are all simply connected.
\end{proof}

\noindent 
Now
we study the components $\mathcal M$ of the moduli parameterizing the surfaces of general type that appear in Theorems~\ref{theorem.construct.surfaces} and~\ref{2to1deformsto2to1}. We compute their dimension and the dimension of the ``hyperelliptic locus'' of $\mathcal M$ for the case of the surfaces in Theorem~\ref{theorem.construct.surfaces}.
We start with some results stated in a more general setting:

\begin{lemma}\label{moduli.dim.lemma}
Let $X$, $Y$ and $\varphi$ be as in~\ref{setup.extra}.
If $\, H^1(\mathcal N_\varphi)=0$ (hence $X$ is unobstructed by Lemma~\ref{varphi.iff.X.unobstructed} and Remark~\ref{canonical.remark}), $Y$ is regular and $p_g(Y)=h^1(\mathcal O_Y(1))=0$, then the (only) component of the moduli containing $[X]$ has dimension \begin{equation}\label{moduli.dim.lemma.equation}
\mu= h^0(\mathcal N_\pi) - h^1(\mathcal N_\pi) + h^1(\mathcal T_Y) - h^0(\mathcal T_Y) + \mathrm{ dim} \;\mathrm{Ext}^1 (\Omega_Y,\omega_Y(-1)).
\end{equation}
\end{lemma}

\begin{proof}
Recall the following sequence of~\cite[(3.3.2)]{Gon}:
\begin{equation}\label{Miguel.seq}
 0 \longrightarrow \mathcal N_\pi \longrightarrow \mathcal N_\varphi \longrightarrow \pi^*\mathcal I/\mathcal I^2 \longrightarrow 0.
\end{equation}
Pushing~\eqref{Miguel.seq} to $Y$ and taking global sections gives
\begin{equation*}
 h^0(\mathcal N_\varphi)=h^0(\mathcal N_\pi) + \textrm{ dim Hom } (\mathcal I/\mathcal I^2,\mathcal O_Y) + \textrm{ dim Hom }(\mathcal I/\mathcal I^2,\omega_Y(-1)) - h^1(\mathcal N_\pi).
\end{equation*}
Since $p_g(Y)=q(Y)=h^1(\mathcal O_Y(1))=0$ (see~\eqref{vanish.for.isom.Hom.Ext}) and $Y$ is embedded in $\mathbf P^N$ by a complete linear series, Lemma~\ref{isom.Hom.Ext} implies that Hom$(\mathcal I/\mathcal I^2,\omega_Y(-1))$ and Ext$^1 (\Omega_Y,\omega_Y(-1))$ are isomorphic so in the above formula we can write dimension of Ext$^1 (\Omega_Y,\omega_Y(-1))$ instead of dimension of Hom$(\mathcal I/\mathcal I^2,\omega_Y(-1))$. On the other hand, since $H^1(\mathcal N_\varphi)=0$, $\varphi$ is unobstructed, so the base of the universal deformation space of $\varphi$ has dimension $h^0(\mathcal N_\varphi)$. Then, the dimension
at $[X]$ of the moduli is $h^0(\mathcal N_\varphi) -$dim PGL$(\mathbf P^N)$. On the other hand, $\textrm{ dim Hom } (\mathcal I/\mathcal I^2,\mathcal O_Y)=h^0(\mathcal N_{Y/\mathbf P^N})=h^0(\mathcal T_{\mathbf P^N|Y})+h^1(\mathcal T_Y)-h^0(\mathcal T_Y)$, because $h^1(\mathcal T_{\mathbf P^N|Y})=0$, what can be easily checked by using the Euler sequence on $\mathbf P^N$ and having in account that $p_g(Y)=h^1(\mathcal O_Y(1))=0$. Finally, since $Y$ is regular, the Euler sequence on $\mathbf P^N$ restricted to $Y$ implies that $h^0(\mathcal T_{\mathbf P^N|Y})=(N+1)^2-1$, which is the dimension of PGL$(\mathbf P^N)$.
\end{proof}

\begin{corollary}\label{hyperelliptic.locus.dim.codim.cor}
Let $X$ be as in Lemma~\ref{moduli.dim.lemma}. Assume in addition that $h^1(\mathcal N_\pi)=0$ and that $Y$ is unobstructed in $\mathbf P^N$. Let $\mathcal M$ be the moduli component of $X$. Then the only irreducible component through $[X]$ of the stratum $\mathcal M_2$ of $\mathcal M$ parameterizing surfaces whose canonical map is a degree $2$ morphism has dimension
\begin{equation*}
 \mu_2=h^0(\mathcal N_\pi) - h^1(\mathcal N_\pi)
  + h^1(\mathcal T_Y) - h^0(\mathcal T_Y)
\end{equation*}
and codimension
\begin{equation*}
 \mu-\mu_2=\mathrm{ dim \, Ext}^1 (\Omega_Y,\omega_Y(-1)).
\end{equation*}
\end{corollary}

\begin{proof}
Recall that $q(Y)=p_g(Y)=0$, $h^1(\mathcal O_Y(1))=0$  by~\eqref{vanish.for.isom.Hom.Ext} and $Y$ is unobstructed by assumption; therefore the hypotheses (1), (2), (3) and (5) of Theorem~\ref{Psi2=0} are satisfied. 
Moreover, hypothesis (4) of Theorem~\ref{Psi2=0} 
follows from Lemma~\ref{normal.pi}, ~\eqref{proving.h1normal.zero} and the assumptions $q(Y)=h^1(\mathcal N_\pi)=0$. Therefore, from the proof of Theorem~\ref{Psi2=0} (see Remark~\ref{kernel.psi2}) it follows that the base of the universal deformation space of $\varphi$ has a stratum parameterizing pairs $(X',\psi)$, where $X'$ are surfaces whose canonical map $\psi$ is a degree $2$ morphism. Furthermore, this stratum is smooth at $[X,\varphi]$ and the dimension of its tangent space is $h^0(\mathcal N_\varphi)-\textrm{ dim Hom } (\mathcal I/\mathcal I^2,\omega_Y(-1))$ (see again Remark~\ref{kernel.psi2}). Letting PGL$(\mathbf P^N)$ act we get a stratum of $\mathcal M$ parameterizing surfaces whose canonical map is a  morphism of degree $2$, whose (only) irreducible component $\mathcal M_2$ passing through $[X]$ has as codimension the dimension of $\textrm{Ext}^1 (\Omega_Y,\omega_Y(-1))$. Then Lemma~\ref{moduli.dim.lemma} implies that the dimension of $\mathcal M_2$ is $h^0(\mathcal N_\pi) - h^1(\mathcal N_\pi)
 + h^1(\mathcal T_Y) - h^0(\mathcal T_Y)$.
\end{proof}

\noindent Now we use Lemma~\ref{moduli.dim.lemma} and Corollary~\ref{hyperelliptic.locus.dim.codim.cor} to compute the dimension of the components of the moduli parameterizing the surfaces of general type of Theorems~~\ref{theorem.construct.surfaces} and~\ref{2to1deformsto2to1}:

\begin{proposition}\label{prop.moduli.2to1deformsto2to1}
 Let $Y$ be $\mathbf P^2$ blown--up at $s$ points in general position and embedded by $|dH-E_1-\cdots-E_s|$ (see~\ref{notation.blowups}) and let $X$ be a surface of general type as in Theorem~\ref{2to1deformsto2to1}.  Then there is only one irreducible component of the moduli containing $[X]$ and its dimension is
\begin{equation*}
 \mu=2d^2+15d+19-8s.
\end{equation*}
\end{proposition}

\begin{proof}
We compute the right--hand--side of~\eqref{moduli.dim.lemma.equation}.
First, recall that Corollary~\ref{exist.abstr.carpets} (1) says that
\begin{equation}\label{Ext1.vanish}
 \textrm{Ext}^1 (\Omega_Y,\omega_Y(-1))=0.
\end{equation}
 Second, recall that in the proof of Theorem~\ref{2to1deformsto2to1} we showed the  vanishing of $H^1(\omega_Y^{-2}(2))$; then~\eqref{proving.h1normal.zero} and the fact that $p_g(Y)=0$ imply
\begin{equation}\label{H1.normal.pi.vanish}
 H^1(\mathcal N_\pi)=0.
\end{equation}
Third, by the same argument used in the proof of Proposition~\ref{vanish.normal}, $h^2(\mathcal T_Y)=0$ (see~\eqref{h2.tangent.Y=0}), so
\begin{equation}\label{Xi.TY}
 h^0(\mathcal T_Y)-h^1(\mathcal T_Y)=\chi(\mathcal T_Y)=2K_Y^2-10\chi(\mathcal O_Y)=8-2s.
\end{equation}

\smallskip
\noindent Finally, to complete the computation of $\mu$ we find the value of $h^0(\mathcal N_\pi)$.
For this, we recall that, by~\eqref{cohom.normalpi.formula},
\begin{equation}\label{h^0.normal.pi}
 H^0(\mathcal N_\pi) \simeq H^0(\mathcal O_B(B)).
\end{equation}
Now, to compute $H^0(\mathcal O_B(B))$ we consider the sequence
\begin{equation}\label{OBB.sequence}
 0 \longrightarrow H^0(\mathcal O_Y) \longrightarrow H^0(\mathcal O_Y(B)) \longrightarrow H^0(\mathcal O_B(B))\longrightarrow 0.
\end{equation}
Arguing as in the proof of Proposition \ref{vanish.normal}, we see that $H^0(\mathcal O_Y(B))=H^0(\mathcal O_{\mathbf P^2}(2d+6) \otimes \mathfrak m^4 )$ is the kernel of the map $\beta$ of~\eqref{beta.restriction.to.4S}.
In the proof of Theorem~\ref{2to1deformsto2to1} we checked that $H^1(\mathcal O_{\mathbf P^2}(2d+6) \otimes \mathfrak m^4)$ vanishes, so $\beta$ is surjective and $h^0(\mathcal O_{\mathbf P^2}(2d+6) \otimes \mathfrak m^4)=\frac{(2d+8)(2d+7)}{2}-10s$. Therefore using this, \eqref{h^0.normal.pi} and~\eqref{OBB.sequence} we get
\begin{equation}\label{h0.normal.pi.dim}
 h^0(\mathcal N_\pi)=h^0(\mathcal O_B(B))=h^0(\mathcal O_Y(B))-1=h^0(\mathcal O_{\mathbf P^2}(2d+6) \otimes \mathfrak m^4)-1=2d^2+15d+27-10.
\end{equation}
Then plugging~\eqref{Ext1.vanish}, \eqref{H1.normal.pi.vanish}, \eqref{Xi.TY} and~\eqref{h0.normal.pi.dim} in~\eqref{moduli.dim.lemma.equation} yields the result.
\end{proof}

\begin{proposition}\label{prop.moduli.2to1deformsto1to1}
 Let $Y$ be $\mathbf P^2$ blown--up at $s$ points in general position and embedded by $|dH-E_1-\cdots-E_s|$ (see~\ref{notation.blowups}) and let $X$ be a surface of general type as in Theorem~\ref{theorem.construct.surfaces}.  Then there is only one irreducible component $\mathcal M$ of the moduli containing $[X]$  and there is only one irreducible component through $[X]$ of the stratum $\mathcal M_2$ of $\mathcal M$ parameterizing surfaces whose canonical map is a degree $2$ morphism. The dimension of $\mathcal M$ is
\begin{equation*}
 \mu=d^2+15d+20-6s
\end{equation*}
 and the dimension of $\mathcal M_2$ is
\begin{equation*}
\mu_2=2d^2+15d+19-8s.
\end{equation*}
 We give $\mu$ and $\mu_2$ explicitly for each value of $d$ and $s$ in the following table:
\bigskip

\centerline{\vbox{\tabskip=0pt \offinterlineskip
%\tabskip= .25 truecm
\def\tablerule{\noalign{\hrule}}
\halign to 3.65truecm
%{\valign to125pt}
{\strut
#& \vrule#
\tabskip=0em plus 3em
%\tabskip= .25 truecm
& \hskip .25cm
\hfil #
\hfil  \hskip .05cm
& \vrule #
& \hskip .3cm
\hfil #
\hfil  \hskip .05cm
& \vrule \vrule \vrule #
&\hskip .3cm
\hfil# \hskip .05cm
&  \vrule #
& \hskip .25cm
\hfil #
\hfil  \hskip .05cm
& \vrule #&
\hskip .15cm \hfil# %\hskip .05truecm
\hfil & \vrule#\tabskip=0pt
\cr\tablerule
%&&Type
&&%\omit\hidewidth
$d$ %\hidewidth
&& %\omit\hidewidth
$s$ %\hidewidth
&& %\omit\hidewidth
$\mu$ %\hidewidth
&&%\omit\hidewidth
$\mu_2$ %\hidewidth
&\cr\tablerule
%\cr
\tablerule
\tablerule
&&$3$
&&$5$
&&$44$
&&$42$
&\cr
%&\vskip -2truecm&\vskip -2truecm&\vskip -2truecm&\vskip -2truecm&\vskip -2truecm&\vskip -2truecm&\vskip -2truecm&\vskip -2truecm&\vskip -2truecm&\vskip -2truecm&\vskip -2truecm&\vskip -2truecm&\vskip -2truecm&\vskip -2truecm&\cr
\tablerule
%&&2
&&$3$
&&$6$
&&$38$
&&$34$
&\cr
\tablerule
%&&3
&&$4$ %\ *$ \hskip -.2truecm
&&$8$
&&$48$
&&$47$
&\cr
\tablerule
%&&4
&&$4$ % \ *$ \hskip -.2truecm
&&$9$
&&$42$
&&$39$
&\cr
\tablerule
%&&5.1
&&$4$ %\hskip -.15truecm
&&$10$
&&$36$
&&$31$
&\cr
\tablerule
&&$5$ %\hskip -.15truecm
&&$13$
&&$42$
&&$40$
&\cr
\tablerule
%&&6.1
&&$5$ %\hskip -.15truecm
&&$14$
&&$36$
&&$32$
&\cr
\tablerule
\noalign{\smallskip}}}}
\medskip
\end{proposition}

\begin{proof}
Proposition~\ref{vanish.normal} (1) says that $h^1(\mathcal N_\pi)=0$. Then Corollary~\ref{hyperelliptic.locus.dim.codim.cor} tells us that
\begin{equation}\label{first.value.of.mu2}
 \mu_2= h^0(\mathcal N_\pi)
+ h^1(\mathcal T_Y) - h^0(\mathcal T_Y).
\end{equation}
The argument used in the proof of Proposition~\ref{prop.moduli.2to1deformsto2to1} to compute $h^0(\mathcal T_Y) - h^1(\mathcal T_Y)$ holds also here so we again have
\begin{equation}\label{Xi.TY.again}
 h^0(\mathcal T_Y) - h^1(\mathcal T_Y)=8-2s.
\end{equation}
Finally, to compute  $h^0(\mathcal N_\pi)$ we argue as in the proof of Proposition~\ref{prop.moduli.2to1deformsto2to1} and get again
\begin{equation}\label{h0.normal.pi.dim.again}
h^0(\mathcal N_\pi)=2d^2+15d+27-10s.
\end{equation}
Then~\eqref{first.value.of.mu2}, \eqref{Xi.TY.again} and~\eqref{h0.normal.pi.dim.again}  imply
\begin{equation}\label{value.of.mu2}
\mu_2=2d^2+15d+19-8s.
\end{equation}

\noindent To compute $\mu$, according to Corollary~\ref{hyperelliptic.locus.dim.codim.cor}, we only need to add to $\mu_2$ the dimension of \linebreak
$\mathrm{Ext}^1(\Omega_Y,\omega_Y (-1))$. As observed in the proof of Corollary~\ref{exist.abstr.carpets}, the latter is $h^1(p^*\Omega_{\mathbf P^2} \otimes \mathcal O_Y(1))$. Recall now that, since $x_1, \dots, x_s$ are in general position, for the values of $d$ and $s$ assumed in this theorem, $H^1(\mathcal O_{\mathbf P^2}(d) \otimes \mathfrak m)=0$, so  $H^1(p^*\Omega_{\mathbf P^2} \otimes \mathcal O_Y(1))$ can be identified with the cokernel of the map $\alpha$ defined in Lemma~\ref{surjectivity}. Thus it only remains to compute the dimension of the cokernel of $\alpha$ for each value of $d$ and $s$. For this recall that in the proof of Proposition~\ref{vanish.normal} we saw that, under our assumptions on $d$ and $s$, $\alpha$ is injective.
Then the dimension of the cokernel of $\alpha$ is
\begin{equation*}
 \frac{(d+2)(d+1)}{2}-s-3(\frac{(d+1)d}{2}-s)=2s+1-d^2.
\end{equation*}
 This together with~\eqref{value.of.mu2} yields
\begin{equation*}
\mu=d^2+15d+20-6s.
\end{equation*}
\end{proof}

%\medskip

\noindent The following remark, regarding components of the moduli of surfaces of
general type, puts into perspective the new results in this paper by contrasting them with what was known earlier.

\begin{remark}\label{horikawa.remark}
 {\rm Two of the families of surfaces constructed in Theorem~\ref{theorem.construct.surfaces} have invariants $(p_g,c_1^2)=(4,6)$ and $(p_g,c_1^2)=(5,8)$. Surfaces of general type with $(p_g,c_1^2)=(4,6)$ were known to Enriques and Max Noether
and studied in depth by Horikawa in~\cite{Horikawa.small3}, whereas surfaces with $(p_g,c_1^2)=(5,8)$  were studied by Horikawa in Section 5 of~\cite{Horikawa.small4}.

\smallskip
\noindent Indeed, if $d=3$ and $s=6$, then $p_g=4$ and $c_1^2=6$ and, according to the notation of~\cite{Horikawa.small3}, the surfaces parameterized by $\mathcal M_2$ are surfaces of Type Ib and the surfaces parameterized by $\mathcal M \smallsetminus \mathcal M_2$ are surfaces of Type Ia. Then, if $d=3, s=6$, the deformation of a canonical double cover $X$  to a finite  morphism of degree $1$ given in  Theorem~\ref{theorem.construct.surfaces} is the inverse of a specialization like the one described in~\cite{Horikawa.small3} from surfaces of Type Ia to surfaces of Type Ib (see the diagram of page 209 and Theorem 7.2 of~\cite{Horikawa.small3}). Horikawa proved in addition that the moduli number of a surface of Type I is $38$ (see~\cite[Theorems 7.1 and 7.2]{Horikawa.small3}); this is obviously the same number we have computed, by different means, in Proposition~\ref{prop.moduli.2to1deformsto1to1}. Horikawa, however, does not compute the dimension of the stratum parameterizing surfaces of Type Ib, whereas for surfaces of Type Ib we do compute the dimension of $\mathcal M_2$, which turns out to be $34$.

\smallskip
\noindent On the other hand, if $d=3$ and $s=5$, then $p_g=5$ and $c_1^2=8$ and, according to the notation of~\cite[Section 5]{Horikawa.small4}, the surfaces parameterized by $\mathcal M_2$ are surfaces of Type Ib and the surfaces parameterized by $\mathcal M \smallsetminus \mathcal M_2$ are surfaces of Type Ia. Then, if  $d=3, s=5$,  the deformation of a canonical double cover $X$  to an embedding (see Remark~\ref{remarkd3}) given in  Theorem~\ref{theorem.construct.surfaces} is the inverse of a specialization like the one given in~\cite[Theorem 5.1.ii]{Horikawa.small4} from surfaces of Type Ia to surfaces of Type Ib. Horikawa proved in addition that the moduli number of a surface of Type I is $44$ (see~\cite[Theorem 5.1.iv]{Horikawa.small4}); this is obviously the same number we have computed in Proposition~\ref{prop.moduli.2to1deformsto1to1}. Horikawa, however, does not compute the dimension of the stratum parameterizing surfaces of Type Ib, whereas for surfaces of Type Ib we do compute the dimension of $\mathcal M_2$, which turns out to be $42$ in this case.}
\end{remark}

\begin{question}\label{moduli4,8} {\rm Let $X$ be a surface as in Theorem~\ref{theorem.construct.surfaces}, with $d=4$ and $s=8$. Let $\mathcal M_{(7,0,16)}$ be the moduli of surfaces with invariants $p_g=7, q=0, c_1^2=16$. The component $\mathcal M$ of $\mathcal M_{(7,0,16)}$  containing $[X]$ has dimension $48$. A standard argument involving a Hilbert scheme dimension computation yields that the dimension of the component $\mathcal M'$ of $\mathcal M_{(7,0,16)}$ parameterizing $(2,2,2,2)$ complete intersections in $\mathbf P^6$ is also $48$. Are $\mathcal M$ and $\mathcal M'$ the same component of $\mathcal M_{(7,0,16)}$?}
\end{question}

\medskip
\noindent We end this section by noting an interesting phenomenon. We show the existence of  infinitely many moduli spaces which have at least two components, one of them parameterizing surfaces whose canonical map is a degree $2$ morphism and the other of them parameterizing  canonically embedded surfaces. This bears further evidence to the complexity of the moduli of surfaces of general type. First we state the following generalization of~\cite[4.5]{AK} (if $m=4$, then the surfaces appearing in Lemma~\ref{divisor.scroll} are the canonically embedded surfaces obtained by Ashikaga and Konno in~\cite[4.5]{AK}):

\begin{lemma}\label{divisor.scroll}
 Let $a,b,c$ and $r$, $Z=S(a,b,c), H$ and $F$ be as in Notation~\ref{scroll.notation}. Let $m, l \in \mathbf Z$, $m \geq 4$. If $ma+l > 0$ and $(m-3)a+r+l-5 >0$, then the general member of the linear system $|mH+lF|$ on $Z$ is a smooth surface $S$ of general type with very ample canonical bundle. In addition, if $m$ is even, $S$ can  degenerate
to a surface $S'$ whose canonical map is a morphism of degree $2$. The image of the canonical morphism of $S'$ has nonnegative Kodaira dimension if and only if $m \geq 6$ and is a conic bundle if $m=4$.
\end{lemma}

\begin{proof}
 The condition $ma+l > 0$ implies that $mH+lF$ is a very ample divisor whereas \linebreak
$(m-3)a+r+l-5 >0$ implies that $K_Z+S$ is a very ample divisor. Then the first claim follows from adjunction and Bertini's theorem. For the second claim, note that if $m$ is even, the argument in~\cite[4.5]{AK} can be easily adapted  to our situation.
\end{proof}

\begin{theorem}\label{2components} Let $\Xi \subset \mathbf N^2$ be the set consisting of those pairs $(x',y)$ for which there exist surfaces $X$ as in Theorem~\ref{2to1deformsto2to1} and surfaces $S$ as in Lemma~\ref{divisor.scroll} such that $(x',y)=(p_g(X),c_1^2(X))=(p_g(S),c_1^2(S))$.
\begin{enumerate}
\item The set $\Xi$ is infinite.
\item For any $(x',y) \in \Xi$, the moduli space of surfaces with invariants $p_g=x', q=0$ and $c_1^2=y$
has at least two different components $\mathcal M$ and $\mathcal M'$, the former containing $[X]$ and the latter containing $[S]$. The general point of $\mathcal M$ corresponds to a surface whose canonical map is a degree $2$ morphism whereas the general point of $\mathcal M'$ is a surface which can be canonically embedded in projective space.
\item The pairs $(x',y)$ of $\Xi$ lie on the lines of the $(x',y)$--plane  which have equations
 \item[] \begin{equation}\label{implicit2}
 y=6\frac{m-3}{m-2}x'-(m-3)(m+3),
\end{equation}
for any $m \in \mathbf Z$, $m \geq 4$. More precisely, the pairs of $\Xi$ are distributed in the $(x',y)$--plane as follows:
\smallskip

\begin{enumerate}
\item for $m=4$, the pairs of $\Xi$ lying on the line of equation~\eqref{implicit2} are exactly $(x',y)=(9,20), (15,38)$, $(23,62)$ and $(33,92)$;
\item for $5 \leq m \leq 10$, on each line of equation~\eqref{implicit2} there are infinitely many pairs of $\Xi$;
\item for $m \geq 11$, the number of pairs of $\Xi$ lying on each line~\eqref{implicit2} is finite (possibly empty).
\end{enumerate}
\end{enumerate}
\end{theorem}

\begin{proof}
Part (2) is a direct consequence of the definition of $\Xi$, Theorem~\ref{2to1deformsto2to1} and Lemma~\ref{divisor.scroll}. Part (1) follows from (3).
The proof of (3) follows essentially from performing some tedious, elementary arithmetics, so we give just an outline. In the proof of Proposition~\ref{noin3folds} (see~\eqref{implicit}) we proved that if $S$ is as in Lemma~\ref{divisor.scroll}, then $(p_g(S),c_1^2(S))$ lies on a line of equation~\eqref{implicit2}. On the other hand, if $X$ is a surface of Theorem~\ref{2to1deformsto2to1}, $p_g(X)=\frac{1}{2}d^2+\frac{3}{2}d-s+1$ and $c_1^2(X)=2d^2-2s$, with
\begin{equation}\label{sbound}
 s \leq \frac{1}{5}d^2+\frac{13}{10}d+\frac{21}{10} \text{ if
} d \geq 7.
\end{equation}
Plugging  $(p_g(X),c_1^2(X))$ in~\eqref{implicit2} then yields
\begin{equation}\label{sfunctionofd}
s=\frac{m-5}{2(2m-7)}d^2+\frac{9(m-3)}{2(2m-7)}d-\frac{(m-3)^2(m+4)}{2(2m-7)}.
\end{equation}
It is clear that if $m=4$ or $m \geq 11$, there are only finitely many values of $d$ (recall that $d$ is a positive integer) satisfying~\eqref{sbound} and~\eqref{sfunctionofd}, whereas if $5 \leq m \leq 10$, there are infinitely many values of $d$  satisfying~\eqref{sbound} and~\eqref{sfunctionofd}.

\smallskip

\noindent We look now at the case $m=4$. In this case~\eqref{sfunctionofd} becomes  $s=-\frac{1}{2}d^2+\frac{9}{2}d-4$, so  $p_g(X)=d^2-3d+5$ and if $X$ is a surface of Theorem~\ref{2to1deformsto2to1} such that $(p_g(X),c_1^2(X))$ lies on the line~\eqref{implicit2}, then $c_1^2(X)=3d^2-9d+8$. Moreover, $s \geq 1$, so we have $d^2-9d+10 \leq 0$ and therefore $2 \leq d \leq 7$. If $d=2$, Theorem~\ref{2to1deformsto2to1}  tells us that $s=1$, so we may rule out this case. If $d=3$, then $s=5$ and this is not allowed by Theorem~\ref{2to1deformsto2to1}, so we may rule out this case also. Then the only possible values of $d$ left are $d=4,5,6$ and $7$. For these, $(d,s)=(4,6), (5,6), (6,5)$ and $(7,3)$ and $(p_g(X),c_1^2(X))=(9,20), (15,38), (23,62)$ and $(33,92)$. Now we see the existence of $S$ as in Lemma~\ref{divisor.scroll}  such that $(p_g(S),c_1^2(S))=(9,20), (15,38), (23,62)$ or $(33,92)$. For $m=4$ it is clear that there exist integers $r$ and $l$ that solve~\eqref{parametric} for $(p_g(S),c_1^2(S))=(9,20), (15,38), (23,62)$ or $(33,92)$. Moreover, by checking case by case, one can find values of $r$ and $l$ and integers
$a, b$ and $c$ as in Notation~\ref{scroll.notation} such that $4a+l > 0$ and $a+r+l-5 >0$ (for instance, for $p_g(S)=9$ we may choose $r=-l=24$ and $a=b=c=7$). Then Lemma~\ref{divisor.scroll} shows the existence of smooth surfaces $S$ with very ample canonical divisor such that
$(p_g(S),c_1^2(S))=(9,20), (15,38), (23,62)$ and $(33,92)$.

\smallskip

\noindent Now fix $5 \leq m \leq 10$.  Recall that if $S$ is a surface as in Lemma~\ref{divisor.scroll}, then $(p_g(S),c_1^2(S))$ satisfy~\eqref{parametric}. Note that the equation
\begin{equation}\label{kappa}
\kappa=\frac 1 6 (m-2)(m-1)(rm+3l)-\frac 1 2 (m-2)(m-1)(m+1)
\end{equation}
has integer solutions $r$ and $l$ for any fixed integer $\kappa$ and, if $\kappa$ is sufficiently large, the equation~\eqref{kappa} has natural solutions $r$ and $l$, $r \geq 6$, and for such solutions we may always find $a, b$ and $c$ as in Notation~\ref{scroll.notation} such that $ma+l > 0$ and $(m-3)a+r+l-5 >0$. Then Lemma~\ref{divisor.scroll} implies that if $x'$ is a sufficiently large integer, then there exist smooth surfaces $S$ with  $(p_g(S),c_1^2(S))=(x',y)$ lying on~\eqref{implicit2}. On the other hand, if $X$ is a surface as in Theorem~\ref{2to1deformsto2to1} and $(p_g(X),c_1^2(X))$ lies on~\eqref{implicit2}, $s$ is a function of $d$ as expressed in~\eqref{sfunctionofd} and
\begin{equation*}
 p_g(X)=\frac{1}{2}d^2+\frac{3}{2}d-s+1.
\end{equation*}
If, in addition, $(p_g(X),c_1^2(X))=(p_g(S),c_1^2(S))$ for some $S$ of Lemma~\ref{divisor.scroll}, then $(p_g(X),c_1^2(X))$ should satisfy~\eqref{parametric}. Putting all this together we conclude that in order for $(x',y)$ to be a point of $\Xi$ it suffices that $(x',y)$ lie on~\eqref{implicit2} and that $x'$ may be written as $\frac{1}{2}d^2+\frac{3}{2}d-s+1$, where $s$ is as in~\eqref{sfunctionofd} and $d$ is a sufficiently large positive integer $d$ satisfying the congruences
\begin{align}\label{congruences}
\text{If } m=5,\ \, &\text{then } d \equiv 1 \text{ or } 2 \ (4) \cr
\text{If } m=6,\ \, &\text{then } d \equiv 0 \text{ or } 3 \ (25) \cr
\text{If } m=7,\ \, &\text{then } d \equiv 4 \text{ or } 6 \ (7) \cr
\text{If } m=8,\ \, &\text{then } d \equiv 5 \text{ or } 19 \ (21) \cr
\text{If } m=9,\ \, &\text{then } d \equiv 6,\  30,\  61 \text{ or } 85 \ (88) \cr
\text{If } m=10, \ \, &\text{then } d \equiv 7,\  20,\  22 \text{ or } 35 \ (39),
\end{align}
(congruences~\eqref{congruences} should be satisfied so that $s$ in~\eqref{sfunctionofd} be a natural number and the system on $r$ and $l$,
\begin{eqnarray*}
x'&=&\frac 1 6 (m-2)(m-1)(rm+3l)-\frac 1 2 (m-2)(m-1)(m+1) \cr
y&=&(m-3)(m-1)(rm+3l)-m(m-3)(3m+1)\cr
\end{eqnarray*}
have integer solutions).
\end{proof}

\begin{remark}
 {\rm If one performs  further computations it is possible to find out, for given $5 \leq m \leq 10$, what the points lying on $\Xi$ and on one of the lines~\eqref{implicit2} exactly are. For instance, if $m=5$, the points $(x',y)$ of $\Xi$ lying on the line~\eqref{implicit2} are precisely those points of~\eqref{implicit2} with $x'=12$ or with $x'$ being an even integer with $x' \geq 16$.

\smallskip

\noindent Also performing some more computations one could find out, for each fixed $m \geq 11$, exactly what points of $\Xi$ lie on the line~\eqref{implicit2}. For instance,
\begin{enumerate}
 \item if $m=11$, then  $(135, 608)$ is the only point that belongs to $\Xi$ and lies on the line~\eqref{implicit2}; it corresponds to the invariants of, on the one hand, surfaces which are general members of the linear system $|11H -7F|$ on $S(1,1,1)$ and, on the other hand, surfaces $X$ in Theorem~\ref{2to1deformsto2to1} with $(d,s)=(20, 96)$;
\item if $m=13$, then $(264, 1280)$ is the only point that belongs to $\Xi$ and lies on the line~\eqref{implicit2}; it corresponds to the invariants of, on the one hand, surfaces which are general members of the linear system $|13H -8F|$ on $S(1,1,1)$ and, on the hand, surfaces $X$ in Theorem~\ref{2to1deformsto2to1} with $(d,s)=(29,201)$;
\item if $m=12, 14, 15, 16, 17$, then there are no points of $\Xi$ lying on the line~\eqref{implicit2}.
\end{enumerate}
}
\end{remark}

\begin{acknowledgement}
{\rm We are very grateful to Edoardo Sernesi for drawing our attention to our earlier work on deformation of morphisms and for suggesting us to use it to construct canonical surfaces. Our grateful thanks also to Madhav Nori, who suggested to us possible applications of our work on deformations to the construction of embedded varieties with given invariants.
We also thank Brian Harbourne for the proof of Lemma 3.3 and for helpful conversations, and Tadashi Ashikaga, for bringing to our attention his result~\cite[4.5]{AK} with Kazuhiro Konno.
We are also very grateful to the referee for his/her comments and corrections, which improved our exposition at various places and made the article more precise.}
\end{acknowledgement}

\end{document}